\documentclass[12pt]{amsart}
\usepackage{amsmath,amssymb,amsthm,enumerate,setspace,graphicx,url,color,multirow,subfig}
\usepackage{natbib}
\usepackage[margin=1in]{geometry}
\theoremstyle{plain}
\newtheorem{theorem}{Theorem}%[section]
\newtheorem{proposition}{\small{Proposition}}%[section]
\theoremstyle{definition}
\newtheorem{lemma}{Lemma}%[section]
\theoremstyle{remark}

\begin{document}

\title{Statistical Tests for Large Tree-structured Data}

\author{Karthik Bharath}
\address{School of Mathematical Sciences, University of Nottingham \newline Nottingham, NG7 2RD, U.K.}
 \email{karthik.bharath@nottingham.ac.uk}

\author{Prabhanjan Kambadur}
\address{Bloomberg LP \newline New York,  NY 10022, USA}
 \email{pkambadur@bloomberg.net}

\author{Dipak. K. Dey}
\address{Department of Statistics, University of Connecticut \newline Storrs, CT 06269, USA}
 \email{dipak.dey@uconn.edu}

\author{Arvind Rao}
\address{Department of Bioinformatics and Computational Biology\newline
 The University of Texas MD Anderson Cancer Center\newline Houston, TX 77030, USA}
 \email{aruppore@mdanderson.org} 
 
 \author{{Veerabhadran Baladandayuthapani}}
\address{Department of Biostatistics, The University of Texas MD Anderson Cancer Center\newline Houston, TX 77030, USA}
 \email{veera@mdanderson.org}
 
\date{}
\maketitle

\begin{abstract}
We develop a general statistical framework for the analysis and inference of large tree-structured data, with a focus on developing asymptotic goodness-of-fit tests. We first propose a consistent statistical model for binary trees, from which we develop a class of invariant tests. Using the model for binary trees, we then construct tests for general trees by using the distributional properties of the Continuum Random Tree, which  arises as the invariant limit for a broad class of models for tree-structured data based on conditioned Galton--Watson processes. 
%The proposed probability models for tree-structured data are easily interpretable and lead to an efficient scheme for simulating a large number of trees with arbitrary sizes. 
The test statistics for the goodness-of-fit tests are simple to compute and are asymptotically distributed as $\chi^2$ and $F$ random variables. We illustrate our methods on an important application of
detecting tumour heterogeneity in brain cancer. We use a novel approach with tree-based representations of magnetic resonance images and employ the developed tests to ascertain tumor heterogeneity between two groups of patients. \\
\\
%\begin{keywords}
 \end{abstract}
%\end{keywords}

\section{Introduction}
Rapid advancements in technology have led to the emergence of datasets in which the underlying quantities of interest are non-Euclidean objects. Increasingly encountered across several disciplines is data that generate tree-like structures or hierarchical representations. 
%These include, for example, database applications that involve the detection of duplicate XML records \citep{YKT, TP}, modeling secondary structures of RNA sequences \citep{SW,LNM, KE}, modeling blood vessels in a human brain \citep{Aydin-Marron, Shen-Marron, WM}, and protein classification \citep{busch}. 
Following a suitable representation of trees, the assumption then is to treat the tree-structured object as an observable quantity that represents the statistical atom. Some central challenges have stymied the systematic development of tools for statistical inference in such settings: the non-Euclidean nature of the set of trees induces an unreasonable dependence of models on the choice of representation and geometry; comparing and sampling trees with different topological and branch length information is not straightforward; the labeling scheme of the vertices and branches influence inference. 

Attempts to address these issues, at least from a modeling perspective, have hitherto been characterized by nonparametric or heuristic/algorithmic approaches specifically motivated by applications. For example, in the context of trees used to model blood vessels in a human brain, \cite{Shen-Marron} used nonparametric functional data analysis methods, \cite{WMA} developed a nonparametric regression model with a tree-structured response variable, and \cite{APWB, AB} used PCA-methods using a suitable metric on trees. \cite{YKT, TP} employed the Edit distance with algorithmic tools for detection of duplicate XML records represented as trees. See also \cite{SW, LNM,KE} for use of algorithmic methods on secondary structures of RNA represented as trees. Despite the flexibility associated with such nonparametric and algorithmic methods for tree-structured data, depending on the choice of representation and metric, developing robust inferential methods, and determining distributions of the statistics of interest are difficult. Critically, without a probability model supporting the simulation of trees through a generative model, it becomes hard to assess the generality of the methods across tree-structured datasets. 

In contrast, parametric probability models for trees were considered by \cite{MS} and \cite{aldous3} wherein the parameters were designed to capture topological features and branch length information. Such probability models form the basis of modeling hierarchical data using Bayesian methods with nonparametric hierarchical priors such as the Dirichlet diffusion tree \citep{RN}, beta-coalescent trees \citep{YH} and stick-breaking priors \citep{AGJ}, and tasks involving clustering \citep{KAH, TDR}. Although sampling classes of trees is convenient with such methods, posterior inference is non-standard, and assessing the quality of posterior samples is difficult. Furthermore, extensions to general classes of trees require substantial revamping of the computational setup; for example, modifying the Dirichlet diffusion tree prior for binary trees to handle trees with an arbitrary number of children requires significant changes in the probability model and computational techniques. 

\begin{figure}[!ht]
\includegraphics[scale=0.5]{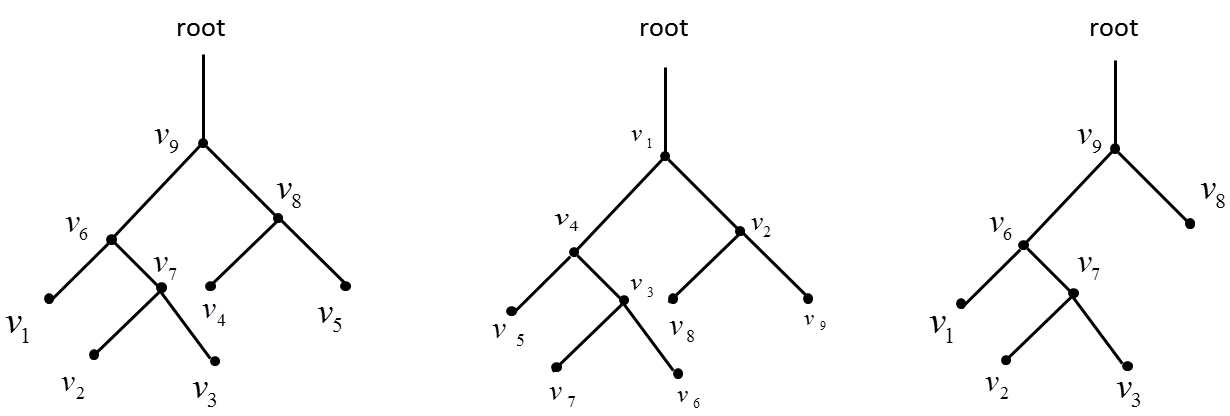}
\caption{\footnotesize Left: Binary tree with $n = 10$ vertices, 9 edges, 5 leaves $(v_1, v_2, v_3, v_4, v_5)$ with a candidate 5-dimensional distribution $f_5$. Middle: Identical binary
	tree with labels of vertices permuted; it is required that the distribution is still $f_5$.
	Right: Resulting binary tree when leaves $v_4$ and $v_5$ are removed from the tree on the left;
	it is required that the distribution becomes the 3-dimensional $f_3$, and is invariant
	to labeling.}
\label{tree_ops}
\end{figure}

Such considerations support an amalgamation of the simplicity of parametric models, in terms of interpretation and simulation, and the flexibility and generality offered by nonparametric or algorithmic methods to develop coherent and principled inferential methods for tree-structured data. This forms the leitmotif of this article. The key challenges in constructing valid probability models for inference are intimately related to the structural components of a tree-structured datum: internal vertices, leaves or terminal vertices and labels. Candidate probability distributions are necessarily multivariate and capture hierarchical relationships between the vertices and the edge length information; in contrast to non-hierarchical relationships, marginalization with respect to any subset of the
vertex set will not be meaningful. Moreover, the distributions should be unaffected by the choice of labelling scheme on the vertices. 
This issue is best understood on binary trees (each node,
except the root, has either 0 or 2 children only) with $n$ leaves, $2n-1$ edges and $2n$ vertices
(including the root and the leaves), when the root has only one child. For simplicity, in this
example let us ignore edge length information. Since the removal of an internal vertex destroys
the tree structure, it is natural to define a distribution using $n$ leaves. As a consequence,
if $f_n$ is an $n$-dimensional distribution on the tree with $n$ leaves, then the removal of one leaf
should ensure that the resulting distribution on the binary tree with $n-1$ leaves is $f_{n-1}$. Additionally, $f_n$ should be \emph{exchangeable} with respect to the leaves, thereby ensuring that the leaf labels are irrelevant. 
Figure \ref{tree_ops} further elucidates these issues using binary trees. 
%A second issue is that when the distribution belongs to a parametric class, the parameters representing distributional properties of the hierarchical relationships ought to retain their
%meaning when leaves are removed or added. Last, since tree-structured objects are not usually observed directly, the labeling scheme used to name vertices, leaves and edges should
%not affect inference; in other words, the probability distribution should be invariant to the
%labelling mechanism. 
Considering the complexity arising from different sources of variation the following requirements are crucial to any statistical method on tree-structured data:\\
	(i) \emph{Uniform sampling}: modeling and inference should be compatible with notions of random sampling from a population of trees;\\
	(ii) \emph{Projective property}: when increasing or decreasing the dimensionality of a tree by
adding or deleting leaves, the probability distribution should be compatible with the notion of obtaining marginal distributions from higher-dimensional ones;\\
	(iiii) \emph{Intepretability of parameters}: parameters should retain interpretability when moving from a subtree to the full tree, or vice versa, when adding or deleting leaves;\\
	 (iv) \emph{Lack of dependence on labeling}: inference should be impervious to labelling of leaves and edges, since the labeling scheme adopted on tree-structured data is usually arbitrary. \\
In view of the requirements, we first propose a generative model for binary trees that satisfies (i)-(iv), resulting in a density that is exchangeable with respect to the leaf labels. Such a property fits into the general framework for exchangeable random structures developed by \cite{OR} in the context of Bayesian models. Using the density, we construct exact parametric one- and two- sample tests that are invariant to permutation of the leaves. This is meaningful because the number of leaves in a binary tree determines its topological structure. The density of binary trees is vital in the generalization to arbitrary trees. We model general trees using the genealogical trees of Galton--Watson branching processes with an offspring distribution, conditioned on the total progeny, and referred to as conditioned Galton--Watson trees (see Section \ref{prelim} for definition). These models are quite versatile and can be used for various types of trees. However, the distribution of a conditioned Galton--Watson tree incorporates only topological information and does not satisfy all requirements (i)-(iv).
To remedy this, we move to the asymptotic limit (as the number of vertices grow without bound) of such trees by seeking recourse to the abstract notion of a limit tree, referred to as the Continuum Random Tree (CRT) from \cite{aldous} and \cite{aldous2}. 
%Analogous to distributions of a continuous-time stochastic process being characterized by the distributions of its finite-dimensional projections, the distribution of the CRT is completely characterized by the density of \emph{randomly} chosen binary subtrees with densities matching the one obtained using the generative model
%for binary trees. 
By suitably modifying and extending the theory of the CRT to accommodate realistic statistical models, we show how the model for binary trees can be extended to a valid probability model for a rich class of general trees. Conditioned Galton--Watson trees, regardless of the offspring distribution, converge weakly to the same limit tree, the CRT; using the two different distributional characterizations of the CRT, we propose goodness-of-fit tests for conditioned Galton--Watson tree models. We design a simulation technique based on an expected linear run-time algorithm by \cite{LD}, which enables the efficient generation of large samples of large conditioned Galton--Watson trees. We then illustrate our methodology through a case study in cancer imaging, where it is well established that a systematic classification of intra- and inter-tumour heterogeneity is crucial for drug development and accurate assessment of response to treatment \citep{NJ}. Obtaining tree-based representations of magnetic resonance (MR) images of brain cancer, we use the proposed tests to detect tumor heterogeneity between cancer patients with long versus short survival times.
 
The paper is structured as follows. In Section \ref{prelim}, we establish the notations and review technical details used in the article. In Section \ref{binary_trees}, we illustrate the motivation for using the CRT through a simple model for binary trees that leads to the definition of a consistent family of densities for the tree-structured data, and the development of invariant tests for tree distributions. In Section \ref{conditioned Galton--Watson}, we review and suitably modify the key ingredients of the CRT, restricting our attention to relevant aspects of their relationships. In Section \ref{GOF}, we construct goodness-of-fit tests for conditioned Galton--Watson tree models using a class of binary subtrees. In view of requirement (iii), we then consider the problem of estimating and interpreting the variance parameter of the offspring distribution of conditioned Galton--Watson tree models. In Section \ref{simulations}, we present results from numerical examples, and implement the tests on a dataset of brain tumor images. In Section \ref{discussion}, we comment on possible extensions and generalizations, and point out some shortcomings of the current framework. 
%%%%%%%%%%%%%%%%%%%%%%%%%%%%%%
\section{Notation and technical preliminaries}
\label{prelim}
In order to make the article self-contained, in this section we detail concepts that figure repeatedly in our methodological development. This section can be skimmed at first reading and referred to later as required. \\
\underline{\noindent\emph{Tree representation}}: A tree is a connected, acyclic graph with a distinguished vertex or node called the root. Nodes are connected through edges that have non-zero lengths. Trees are allowed to have unequal numbers of vertices or nodes. We confine our attention to finite, rooted, ordered trees: trees with a root, containing a finite number of vertices, and having an ordering among the offspring at a parent node. This leads to the notion of a left and right child in the case of binary trees. Ordered or planar representations of unordered trees can be obtained using the scheme of \cite{AldPit}. 
%Strict binary trees are those in which each vertex has 0 or 2 children. To avoid confusion, throughout this article, binary trees indicate strict binary trees. 
Following the notation used by \cite{aldous}, we represent a finite rooted tree  $\tau_n$ with $n$ vertices, including the root, as a point $\tau_n$  in the product space $\mathcal{T}_n \times \mathbb{R}_{+}^{n-1}$, where $\mathcal{T}_n$ is the set of all finite trees on $n$ vertices.  This implies that the connectivity information between the vertices is contained in $\mathcal{T}_n$. For simplicity, we choose the notation
$\tau_n=\Big(\mathcal{V}(\tau_n),\mathcal{E}(\tau_n)\Big)$, where $\mathcal{V}(\tau_n)=(\text{root},v_1,\ldots,v_{n-1})$ is the set of vertices and $\mathcal{E}(\tau_n)=(e_1,\ldots,e_{n-1})$ represents the set of edge lengths. It is tacitly assumed that the notation contains the information that relates each edge to a unique pair of vertices, and  we always view a tree $\tau_n$ as a point in the product space containing topological and edge length information, and use the notation of the vertex and edge length set for convenience. We denote by $\tau_n$ the tree with $n$ vertices including the root, and by $\tau(n)$ the tree with $n$ terminal vertices or leaves.\\
\underline{\emph{Conditioned Galton--Watson trees}}:  
A Galton--Watson process $\{X_n\}_{n \geq 0}$ with offspring distribution $(\pi_k, k=0,1,\ldots)$ is a $\mathbb{Z}_+$-valued discrete-time Markov chain with transition function $P(X_{n+1}=k|X_n=m)=\pi_k^{*m}$, where $*$ denotes convolution operator. If $\mu=\sum_k k\pi_k$, the process $X_n$ is critical if $\mu=1$,
sub-critical if $\mu<1$ and super-critical if $\mu>1$. 
 
\begin{figure}[!ht]
 	\includegraphics[trim=10mm 100mm 15mm 30mm, clip,totalheight=0.2\textheight]{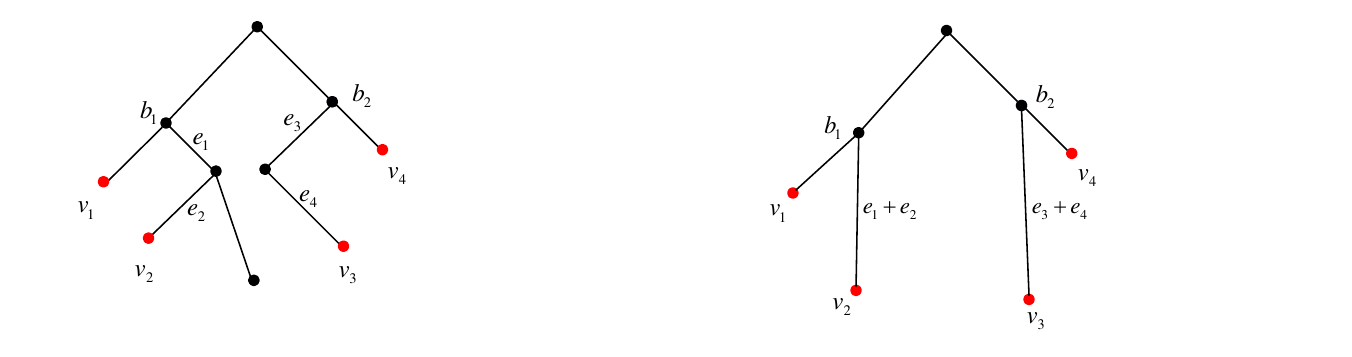}
 	\caption{\footnotesize A tree on the left and its $L$-tree on the right corresponding to vertices $\{v_1,v_2,v_3,v_4\}$. The $L$-tree contains the root, the branch points $b_1$, $b_2$ and the set of vertices $\{v_1,v_2,v_3,v_4\}$. All the vertices chosen are leaves.}
 	\label{LCAtree}
 \end{figure}

A Galton--Watson tree $\tau$ corresponding to $X_n$ is constructed by recursively starting with the root and giving each node $v$, independent of other nodes, a number $o(v,\tau)$ of children with probability $\pi_{o(v,\tau)}$, where $o(v,\tau)\in \{0,1,\ldots\}$ is the out-degree or the number of children of vertex $v$ in the vertex set $\mathcal{V}(\tau)$ of $\tau$.  
%Figure \ref{CGW} provides a graphical description. 
%\begin{figure}[!ht]
%\includegraphics[scale=0.5]{figs/CGW}
%\caption{\footnotesize Construction of Galton--Watson trees. (i): A draw from $\pi_k$ gives 1, and 1 vertex is attached to the root; (ii) the second draw gives 2, and 2 vertices are attached as children to $v_1$; (iii) the third draw gives 3, and 3 vertices are attached to $v_2$ as children; (iv) the fourth draw gives 1, and 1 vertex is attached to $v_3$. There is no edge length information in this construction, and no preferred order of attachment.}
%\label{CGW}
%\end{figure}
$\{\pi_k\}$ is referred to as the offspring distribution. When conditioned to have $n$-vertices, the resulting tree $\tau_n$ is known as a conditioned Galton--Watson tree with distribution
%\begin{equation*}
%P(\tau=t)=\prod_{v \in \mathcal{V}(t)}\pi_{o(v,t)},
%\end{equation*}
\begin{equation*}
P(\tau_n=t) \propto \prod_{v \in \mathcal{V}(t)}\pi_{o(v,t)} \quad \text{on  }\{t: \text{ cardinality of } \mathcal{V}(t)=n\}.
\end{equation*}
\\
%\emph{Continuum Random Tree (CRT)}: 
%We give a heuristic definition of the CRT that suffices for our needs and refer the interested reader to \cite{aldous2} for details.  The CRT arises as the `continuous' limit of critical conditioned Galton--Watson trees when the number of vertices is allowed to grow without bound. One way to compare trees of different sizes is to embed $\tau_n=(\mathcal{V}(\tau_n),\mathcal{E}(\tau_n)) $ as an element of the linear space $\ell_1$, the Banach space of infinite sequences $x=(x_1,x_2,\ldots)$ such that $|| x||=\sum_i |x_i| < \infty$. The CRT or the \emph{Brownian} CRT is the unique tree that is the closure (with respect to the norm on $\ell_1)$) of the union $\cup_{n\geq1}\tau_n$.  \\
\underline{\emph{Least Common Ancestor Tree ($L$-tree)}}: 
 For a tree $\tau_n=(\mathcal{V}(\tau_n),\mathcal{E}(\tau_n)) $, we define its least common ancestor tree ($L$-tree) in the following manner: choose a subset $B$ of $\mathcal{V}(\tau_n)$;  for vertices $v_1$ and $v_2$ in $B$ find their last common ancestor, or the branch point after which the paths to $v_1$ and $v_2$ from the root diverge or branch out. Then, the $L$-tree corresponding to the subset $B$ of the vertices of $\tau_n$ is the tree, denoted as $L(\tau_n,B)$, that contains the root, the vertices of B and all the branch points with distances from the root to the vertices of $B$ preserved. Figure \ref{LCAtree} illustrates this with $B=\{v_1,v_2,v_3,v_4\}$; the branch points are $b_1$ and $b_2$, and in order to preserve the distances from the root to the vertices of $B$, the new branch from the branch points to the elements of $B$ has a length that is equal to the sum of the edges along the path from the root to the elements of $B$ in the original tree.
%%%%%%%%%%%%%%%%%%%%%%%%
\section{A generative model and test for binary trees}
\label{binary_trees}
We begin with binary trees for two reasons. We intend to define probability distributions on trees using leaf-level information, which for binary trees determines the topology completely. Second, the distributional properties of the CRT, which we intend to use to define asymptotic models on general trees, are completely characterized by random binary subtrees. Consequently, a systematic development of methods with desirable properties for binary trees provides similarly attractive ones for general trees. The results of this section are pivotal in developing inferential models for general trees. We first describe a generative model and its properties and then construct one- and two-sample tests that are invariant to leaf-labeling. 

Note that a binary tree with $k$ leaves $\tau(k)$ is isomorphic to the binary tree $\tau_{2k}$ with $2k$ vertices. Suppose $t_1,t_2,\ldots$ are the arrival times of a non-homogeneous Poisson process with rate $\lambda(t)=t$. Let $\tau(1)$ be a tree with one leaf formed by attaching an edge of length $t_1$ to the root. This tree will have a density that corresponds to the time until the first event of the non-homogeneous Poisson process, with mean function $m(t)=\int_0^t sds=t^2/2$. A tree with two leaves $\tau(2)$ is then constructed by adding a second edge of length $t_2-t_1$, the inter-arrival time, to a point on $\tau(1)$ chosen according to a uniform distribution on $[0,t_1]$. Recursively, the tree $\tau({k+1})$ is obtained from $\tau(k)$ by attaching an edge of length $t_{k+1}-t_k$ to a point chosen randomly on the tree (i.e., according to a uniform distribution on the sum of the branch lengths of $\tau(k)$) and labeling the new leaf $k+1$. Since we are interested in only ordered trees, an edge can be  attached to the left or to the right of the randomly chosen point with equal probability. Observe that $\tau(k)$ has $2k-1$ branches with lengths, say $(x_1,\ldots,x_{2k-1})$, and $\tau({k+1})$ has $2k+1$ branches of lengths $(y_1,\ldots,y_{2k+1})$, with $\tau({k+1})$ formed by the above construction by splitting a branch of length $x_j$ into two branches $y_{j1}$ and $y_{j2}$ with $x_j=y_{j1}+y_{j2}$. Therefore, through induction, the `density' of $\tau({k+1})$ is 
\begin{align*}
g_{k+1}(&\tau({k+1}))=
g_k(\tau(k)) \times \text{(density of $k+1$th inter-arrival time)} \\
& \times \text{(probability of attaching the new edge to $\tau(k)$)} \times \text{(probability of a left or right attachment)}.
\end{align*}
Note that $g(\cdot)$ at this stage is a density (non-normalized) with respect to the product Lebesgue measure on the edge lengths only, with the topological information of the tree being captured only in the decision to add an edge to the left or right with equal probability. Denoting the sum of the branch lengths of $\tau(k)$ by $s_k$, and noting from the construction that the third term is $1/s_k$ and the fourth is $1/2$, we obtain,
$g_{k+1}(\tau({k+1}))=g_k(\tau(k)) s_{k+1}e^{-\frac{1}{2}(s_{k+1}^2-s_k^2)}\frac{1}{2s_k},$
which reduces to 
\begin{equation*}
g_k(\tau(k))=\frac{1}{2^{k-1}}se^{-s^2/2},\quad s=\displaystyle \sum_{i=1}^{2k-1}x_i.
\end{equation*}
 The function $g(\cdot)$ is non-negative but does not integrate to one. By virtue of the construction mechanism, since at each step $i$ there are $2i-1$ possible tree topologies, the structural information of the tree is captured in the normalizing constant:
 $$\int_{x_1}\cdots\int_{x_{2k-1}} se^{-s^2/2}dx_1\ldots dx_{2k-1}=\prod_{i=1}^{k-1}\frac{1}{2i-1}.$$
Therefore, the probability density of the ordered and rooted binary tree with $k$ leaves is
 \begin{equation}\label{nonhomo}
 f_k(\tau(k))=\left[\prod_{i=1}^{k-1}\frac{1}{2i-1}\right]^{-1}g_k(\tau(k))=\left[\prod_{i=1}^{k-1}\frac{1}{2i-1}\right]^{-1}\frac{1}{2^{k-1}}se^{-s^2/2},\quad s=\displaystyle \sum_{i=1}^{2k-1}x_i.
 \end{equation}
 %%%%% FIGURE (BINARY TREE CONSTRUCTION)%%%%%%%%%%
 \begin{figure}[ht]
 	\centering
 	\includegraphics[scale=0.4]{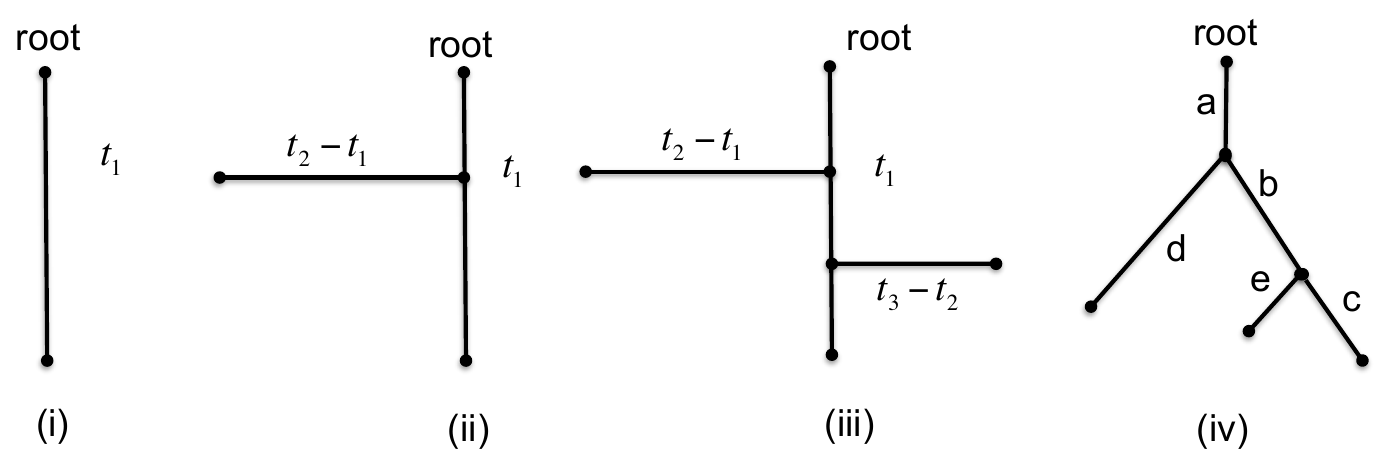}
 	\caption{\footnotesize (i), (ii), (iii): Construction of binary tree $\tau(3)$ with 3 leaves and 2$\times$3-1= 5 edges, 6 vertices (including root) from nonhomogeneous Poisson process model with arrival times $t_1$, $t_2$ and $t_3$. (iv): A binary tree isomorphic to the tree in (iii) where a+b+e=$t_1$, c=$t_3-t_2$ and d=$t_2-t_1$.}
 	\label{binary_poisson}
 \end{figure}
  %%%%%%%%%
 \subsection{Interpretation and properties of the density}
The density $f_k$ can be interpreted as assigning mass to a tree with branch lengths $(x_1,\ldots,x_{2k-1})$ with a canonical binary tree topology on $k$ leaves or $2k$ vertices. That is, the mass is the same for rooted binary trees with $k$ leaves of all possible topologies with branch lengths $(x_1,\ldots,x_{2k-1})$.
We hence observe that, conditional on $k$, $f_k$ is a density with respect to the product measure $u_k \otimes dx$ on $\mathcal{T}_{2k} \times \mathbb{R}_{+}^{2k-1}$, where $u_k$ is a uniform measure on all rooted binary trees on $k$ leaves, and $dx$ is the Lebesgue measure on $\mathbb{R}_{+}^{2k-1}$. The function $f_k$ is interpreted as a density conditional on $K=k$, the number of inter-event times of the Poisson process considered. Indeed, $K$ is Poisson distributed with mean $s^2/2$.

The model is particularly useful in applications where branch lengths are perhaps more important than topology \citep{FL, RL}.  We introduce a parametric form of the density obtained by defining $\lambda(t,\theta)=\theta t$ that leads to a parametric family $\{f_{k,\theta}; \theta \in \Theta\}$, wherein trees of varying branch lengths can be generated by varying $\theta $ in $\Theta$. The reason for the simple multiplicative form lies in our insistence that the density satisfy properties (i)-(iv) in the Introduction. The projective property of the induced probability measures (via the construction) is easily violated when complicated parametric forms are considered for the rate function of the Poisson process. Based on the requirements in the Introduction, we briefly summarize four attractive properties of the density $f_{k,\theta}$ as a statistical model. Formal statements of the properties and their proofs are available in the Section 2.1 (Theorem 2.1 and Proposition 2.1) of Supplementary Material. \\
(i) The uniform measure $u_k$ on the topology ensures that the particular binary tree has been randomly sampled from a population of binary trees on $k$ leaves. \\
(ii) For every $\theta \in \Theta$ $f_{k,\theta}$  is compatible with marginalization over the number of leaves $k$. \\
(iii) For $\tau(k)$, the density $f_{k, \theta}$ depends on the branch lengths $x_1,\ldots,x_{2k-1}$ only through their sum $\sum x_i$, which is Gamma distributed with shape $k$ and scale 2. \\
(iv) From (\ref{nonhomo}), it is evident that $f_{k,\cdot}$ is exchangeable with respect to the branch lengths. 
\subsection{Invariant One- and Two-sample tests}
Hypothesis tests for trees entail testing for a generative probability model in which the probability distribution depends on the number of vertices and branches. When we state that an independent sample is generated from a density, we are essentially viewing the trees as independent trees generated from the same probability model, regardless of the number of vertices. Since the density $f_{k,\theta}$ is uniform on all topologies of $k$-leaved binary trees, for a chosen tree topology, exchangeability of the leaves implies that the data space $\mathcal{X}_k:=\mathbb{R}_{+}^{2k-1}$ is invariant to the symmetric group $G$ of \emph{all} possible permutations of the branch lengths. Thus every $g \in G$ acting on $\mathcal{X}_k$ induces an action $g'$ on $\Theta$; that is, if $f_{k,\theta'}$ is the resulting density following the transformation $g\mathcal{X}_k$ (with density $f_{k,\theta}$), then $g'\theta=\theta'$. It is easy to observe from the Poisson process construction that $g'$ is the identity function (i.e., $g'\theta=\theta'=\theta$ for every $\theta \in \Theta$), and therefore $\Theta$ is also invariant under any $g \in G$. This leads to a testing procedure for binary trees that is invariant to the action of $G$ on the leaves. Technical details ensuring this are in Proposition 2.2 of Supplementary Material. 
%\begin{remark}
%There is a difference between the total path length that characterizes $f$ and the $k$th arrival time that characterizes the process. Recall that while constructing $\tau(k+1)$ from $\tau(k)$, the new edge with length given by the inter-arrival time is attached uniformly on $\tau(k)$ with total path length $s$. As a consequence, once the number of leaves is observed, the density of the tree arising from the construction is fully characterized by the total path length. Generally speaking, in contrast to the homogeneous Poisson process with a constant rate function, the distribution of the $k$th arrival time, in general, does not characterize a non-homogeneous Poisson process. However, it is always possible to convert a non-homogeneous Poisson process to a homogeneous one with a unit rate by warping the time scale: in our setup, if $N(t)$ is the non-homogeneous Poisson process with intensity function $m(t)=t^2/2$ (rate $\lambda(t)=t$), then the process $N(m^{-1}(t))=N(\sqrt{2t})$ is a homogeneous Poisson process with a rate of one. Under the non-linear time scale of the homogeneous Poisson process, in order to obtain a density that depends solely on the total path length, we need to appropriately modify the step in the construction that involves choosing a point on $\tau(k)$ at which to attach a new edge to obtain $\tau(k+1)$.
%\end{remark}
\begin{theorem}
	\label{gof}
	 Suppose $\boldsymbol{\tau(n)}=(\tau({n_1}),\ldots,\tau({n_p}))$ and $\boldsymbol{\eta(m)}=(\eta({m_1}),\ldots,\eta({m_q}))$ are independent samples of binary trees from $f_{\theta_0}$ and $f_{\theta_1}$ respectively. 
	\begin{enumerate}
		\item Consider the critical function
		\begin{equation*}
		\phi(\theta',\mathbf{n},\alpha)=\left\{
		\begin{array}{l l}
		1 & \quad \text{if  }  \displaystyle\theta'\sum_{i=1}^p s_i <\chi_{\alpha,2\sum_{i=1}^pn_i} \quad or \quad 
		 \displaystyle\theta'\sum_{i=1}^p s_i >\chi_{1-\alpha,2\sum_{i=1}^pn_i};  \\
		0 & \quad \text{otherwise} ,
		\end{array}
		\right .
		\end{equation*}
		where $s_i$ is the sum of the branch lengths of $\tau({n_i})$ and $\chi_{\alpha,b}$ denotes the $\alpha$th percentile of a $\chi^2$ distribution with $b$ degrees of freedom. For the hypotheses $H_0: \theta_0=\theta'$ against $H_1:\theta_0 \neq \theta'$, $\phi(\cdot)$ is invariant to the permutation of leaves, and is such that $E_{H_0}\phi(\theta',\mathbf{n},\alpha)=\alpha$.
		\item Let $r_j$ denote the sum of the branch lengths of $\eta({m_j})$. Then, the critical function
		\begin{equation*}
		\psi(\mathbf{n,m},\alpha)=\left\{
		\begin{array}{l l}
		1 & \quad \text{if  }  \frac{\sum_{i=1}^p s_i }{\sum_{j=1}^q r_j }>\left(\frac{\sum_{i=1}^pn_i}{\sum_{j=1}^qm_j}\right)F_{1-\alpha,2\sum_{i=1}^pn_i,2\sum_{j=1}^qm_j};\\
		0 & \quad \text{otherwise} ,
		\end{array}
		\right .
		\end{equation*}
		where $F_{\alpha,a,b}$ is the $\alpha$th percentile of an $F$ distribution with $a$ and $b$ degrees of freedom, for testing $H_0:\theta_0=\theta_1$ against $H_1:\theta_0\neq \theta_1$, and is such that $E_{H_0}\psi(\mathbf{n},\mathbf{m},\alpha)=\alpha$.
	\end{enumerate}
\end{theorem}
In Theorem \ref{gof} for the particular case in which $\theta'=1$ without explicitly specifying an alternative hypothesis, the test is a goodness-of-fit test for the non-homogeneous Poisson process model generating trees with intensity $\lambda(t)=t$. It can easily be checked that the invariance to the action of the symmetry group on the leaves is satisfied.

%%%%%%%%%%%%%%%%%%%%%%%%%%%%%%%%%%%
\section{Conditioned Galton--Watson tree models and the CRT}
\label{conditioned Galton--Watson}
The tests in Theorem \ref{gof} can neither be used to detect topological variation within or across samples, nor to test non-binary trees. In this section, we describe the ingredients that allow us to extend the use of the density in (\ref{nonhomo}) to non-binary trees.
\subsection{Conditioned Galton--Watson tree models}
Conditioned Galton--Watson trees are genealogical trees of Galton--Watson processes conditioned on the total progeny. We refer to the definition in Section \ref{prelim} of a conditioned Galton--Watson tree $\tau_n$ with $n$ vertices. \emph{The key property of interest is that for a fixed offspring distribution $\pi_k$, the corresponding conditioned Galton--Watson tree can be viewed as being picked according to a uniform distribution on certain types of trees with $n$ vertices}. For example, if we wish to choose a binary tree with $n$ vertices according to a uniform distribution on the space of $n$-vertex binary trees, we can equivalently construct a conditioned Galton--Watson tree with an offspring distribution of $1-p$ and $p$ each for 0 and 2 children, respectively, where $p$ is the probability of having $2$ children. 
Other examples include geometric offspring distribution with probability 1/2 for ordered trees with unrestricted degree; Binomial distribution on 2 trails with probability 1/2 for trees with 0,1 or 2 vertex degrees; uniform distribution on $\{0,1,2\}$ for unordered and unlabelled trees with 0,1 or 2 vertex degrees . For a detailed look at the properties of conditioned Galton--Watson trees, and their relationship with other tree models, we refer the interested reader to \cite{janson}.

%
%Conditioned Galton--Watson trees comprise a large class of trees that offer flexibility from a modeling perspective. We enumerate a few that offer a broad range of models for tree-structured data based on specific parameterizations.
%\begin{enumerate}[(i)]\label{sigma}
%\item Ordered trees with unrestricted degree: conditioned Galton--Watson trees with the geometric offspring distribution with success probability $1/2$ (Geo(1/2).
%\item Binary trees: conditioned Galton--Watson trees that contain 0,1 or 2 children according to a binomial distribution, with 2 trials and success probability 1/2 (Bin(2,p)).
%\item Strict binary trees that are ordered: conditioned Galton--Watson trees where the vertices contain either 0 or 2 children with equal probability (Bin(p)).
%\item Unary-binary trees that are ordered: conditioned Galton--Watson trees with vertices that contain 0, 1 or 2 children with equal probability (UB(1/3,1/3)).
%\item Unary-binary trees that are unordered and unlabelled:  conditioned Galton--Watson trees with vertices containing 0, 1 or 2 children with respective probabilities $\pi_0=\frac{1}{2+\sqrt{2}}$, $\pi_1=\frac{\sqrt{2}}{2+\sqrt{2}}$ and $\pi_3=\frac{1}{2+\sqrt{2}}$.
%\item $m$-ary trees: conditioned Galton--Watson trees with binomial offspring distribution with $m$ trials and success probability $1/m$. 
%\end{enumerate}
Evidently, conditioned Galton--Watson trees can be constructed from critical, sub-critical and super-critical Galton--Watson process. However, asymptotic behavior of such trees (as number of vertices tends to infinity) blurs the division between between trees from critical and non-critical Galton--Watson processes: 
in the sub-critical or super-critical case with offspring distribution $\pi_i$ , as long as there exists a $\lambda>0$ such that $\sum_{i\geq0}\pi_i\lambda^i <\infty$, the asymptotic behavior resembles that of a critical Galton--Watson process with finite offspring variance. For modeling purposes such a property is useful since a broad class of distributions can be used as offspring distributions for the Galton--Watson process, giving rise to a rich class of trees.  The following Proposition clarifies this for a few offspring distributions. 
%Despite their general nature, two issues are immediately apparent: there is no branch length information in the conditioned Galton--Watson tree models; offspring distributions have unit mean (critical Galton--Watson process). As potential models for real tree-structured data, we require more. %For example, it is not enough that the only potential model for strict binary trees is the conditioned Galton--Watson tree with offspring distribution 1/2 for both 0 and 2 children; similar arguments apply to the other cases. 
%The following proposition demonstrates the possibility of extending the class of potential offspring distributions to certain sub- and super-critical cases for modeling. 
\begin{proposition}\label{binaryp}
Conditioned Galton--Watson trees with the following offspring distributions can be modeled as critical conditioned Galton--Watson trees for every $0<p<1$ in (i)-(iii) and $0<p_0,p_1<1$ in (iv).
\begin{enumerate}[(i)]
\item $\pi_i=(1-p)^{i-1}p$ for  $i=1,2,\ldots$ and $0<p<1$ ;
\item $\pi_i=\frac{2!}{(2-i)!i!} p^i(1-p)^{2-i}$ for $i=0,1,2$ and $0<p<1$;
\item $\pi_0=1-p, \pi_2=p$ and $0<p<1$;
\item $\pi_0=p_0, \pi_1=p_1, \pi_2=1-p_0-p_1$ and $0<p_0,p_1<1$.
\end{enumerate}
\end{proposition}
%\noindent Proposition \ref{binaryp} clarifies the assumption of criticality for the convergence of conditioned Galton--Watson trees to the CRT: in the sub-critical or super-critical case with offspring distribution $\pi_i$ , as long as there exists a $\lambda>0$ such that $\sum_{i\geq0}\pi_i\lambda^i <\infty$, the asymptotic behaviour resembles that of a critical Galton--Watson process with finite offspring variance.
 \noindent At this point we may choose to view inference on conditioned Galton--Watson trees as inference on Galton--Watson processes; a good source for available methods is \cite{PG}. However, we move to the asymptotic setting for at least two reasons: to incorporate branch length information within the conditioned Galton--Watson tree models; and to obtain knowledge of the distributions of local structural aspects like height and variations in branching structure, through weak convergence techniques based on the CRT.

\subsection{Continuum Random Tree and $L$-tree models}
 A terse and heuristic definition of the CRT is given in Section \ref{prelim}; see \cite{aldous2} for a formal definition. The CRT arises as a `continuous' limit of conditioned Galton--Watson trees as the number of vertices tend to infinity, regardless of the offspring distribution (upto a scaling factor). In this setting, the variance parameter $\sigma^2$ of the offspring distribution of conditioned Galton--Watson trees appears in the limit. The distribution of the CRT can be characterized in two equivalent ways: as limit of randomly chosen binary subtrees of conditioned Galton--Watson trees; as a weak limit of a continuous function constructed by a walk on conditioned Galton--Watson trees. Both characterizations will be profitably used to construct goodness-of-fit tests. 

The class of binary subtrees that characterize the distribution of the CRT  arises as the limit of randomly chosen binary subtrees of conditioned Galton--Watson trees known as $L$-trees used in various applications \citep{GM, ASSU}; for a constructive definition of an $L$-tree see Section \ref{prelim}. For a conditioned Galton--Watson tree $\tau_n$, fix $k<n$ and choose $k$ \emph{leaves} according to a uniform distribution on $\mathcal{V}(\tau_n)$.  %For our purposes, we randomly choose \emph{only} the leaves or the terminal nodes from $\mathcal{V}(\tau_n)$ since $L$-trees constructed in this manner have a direct relationship with the CRT. 
Therefore, the $L$-tree of a conditioned Galton--Watson tree that is a binary tree with $k$ leaves and $2k-1$ edges. The  following characterization of the CRT based on the limit of the $L$-trees, modified suitably for the case of ordered conditioned Galton-Watson trees, is useful for our purpose. The limit binary trees can be viewed as ``marginals" of the CRT or its finite-dimensional projections. 
\begin{lemma}\citep{aldous2})
	\label{marginal}
From a conditioned Galton--Watson tree $\tau_n=(\mathcal{V}(\tau_n),\mathcal{E}(\tau_n))$ generated from an offspring distribution with variance $\sigma^2 <\infty$, for a fixed $k<n$, let $L(\tau_n,\{l_1,\ldots,l_k\})$ be an $L$-tree obtained by leaves $\l_i \in \mathcal{V}(\tau_n)$ with  $s=e_1+\cdots+e_{2k-1}$. Then, as $n \to \infty$, for a fixed $k$, there exists a consistent family $(\mathcal{C}(k), k \geq 1)$ of binary trees with $k$ leaves that define the CRT that has the density
\begin{equation}\label{LCAdensity}
f_{k,\sigma^2}(c(k))=\left[\prod_{i=1}^{k-1}\frac{1}{2i-1}\right]^{-1}\frac{1}{2^{k-1}}(\sigma^2)^k s e^{\frac{-s^2\sigma^2}{2}}\qquad \sigma^ 2 \in \mathcal{S}.
\end{equation}
\end{lemma}
The $L$-trees provide the link through which the density for binary trees obtained using the Poisson process construction is related to conditioned Galton--Watson trees and the CRT. The density in (\ref{LCAdensity}) coincides with density obtained for binary trees in (\ref{nonhomo}) from the non-homogeneous Poisson process construction, and can be used to approximate the distribution of $L$-trees of conditioned Galton--Watson trees from \emph{any} offspring distribution with a finite variance $\sigma^2$. Bearing in mind that the binary tree model satisfies requirements (i)-(iv), the idea then is to model $L$-trees of observed trees with the probability density in (\ref{LCAdensity}), conduct inference, and then extend the resulting conclusions to the whole tree. 
 
The construction of an  $L$-tree incorporates both topological and branch length information. Lemma \ref{marginal} implies that the asymptotic distribution of the sequence of $L$-trees of any conditioned Galton--Watson tree model can be approximated by the sequence in (\ref{LCAdensity}). The freedom associated with the choice of the number of leaves used to construct an $L$-tree of a given tree helps in reducing the dimensionality of the inferential problem involving large trees. Two tree populations can be distinguished by constructing low-dimensional summary statistics by choosing a small proportion of the leaves while constructing an $L$-tree, contingent on consistent estimation of $\sigma^2$. An important issue is the interpretability of $\sigma^2$ in the context of using $f_{\cdot, \sigma^2}$ as a parametric statistical model on the $L$-trees. Formal statement of this property and its proof can be found in Proposition 2.3. in the Supplementary Material. 
%\begin{enumerate}
%	\item conditioned Galton--Watson tree models, parameterized by offspring distributions, provide a rich family of models for topologies of tree-structured data of various kinds;
%	\item CRT approximation of these models through binary $L$-trees allow us to incorporate branch length information, in addition to topological information, and leads to the consistent family of densities prescribed in (\ref{LCAdensity}).
%	\item $L$-trees act as dimension-reduction tools for tree-structured data with a very large number of vertices. 
%\end{enumerate}
%At this juncture, we clarify the notations frequently used in the rest of the paper:
%\begin{itemize}
%	\item $\tau_n$: tree with $n$ vertices including the root;
%	\item $\tau(n)$:  tree with $n$ terminal vertices or leaves;
%	\item $L(\tau_n,B)$: an $L$-tree of a tree $\tau_n$ formed by choosing a subset $B$ of the leaves;
%	\item $\sigma^2$: variance parameter of the offspring distribution of conditioned Galton--Watson trees.
%\end{itemize}

\section{Goodness-of-fit tests for conditioned Galton--Watson trees}
\label{GOF}
The employment of conditioned Galton--Watson trees models, through their connection to the CRT, permits us to extend the use of the consistent parametric family of densities on binary trees to more general settings.
%If we wish to construct asymptotic parametric (w.r.t. $\sigma^2$) tests, there is an identifiability issue; for example,  $\sigma^2=2/3$ can imply an offspring distribution with equal probability for 0,1 and 2 children (unary-binary trees), and the binomial distribution with $3$ trials with success probability $1/3$.
Importantly, the CRT, characterized by the limiting $L$-trees, is the invariant limit for conditioned Galton--Watson trees from \emph{any} offspring distribution with finite variance $\sigma^2$. Consequently, the density in (\ref{LCAdensity}) can be used to approximate the density of a binary $L$-tree of \emph{any} conditioned Galton--Watson tree. In contrast to the situation with binary trees, bearing in mind the invariant limit (CRT), the appropriate test would be a goodness-of-fit test for conditioned Galton--Watson trees, where $\sigma^2$ is viewed as a nuisance parameter that needs to be estimated consistently. We therefore develop tests to check whether the data has been generated from a conditioned Galton--Watson tree model. The null hypothesis would be that the samples of trees are independent copies of a conditioned Galton--Watson tree with finite-variance offspring distribution. This generalizes the test in Theorem \ref{gof} which is applicable only to binary trees.  

%Since the limiting density of an $L$-tree coincides with the density obtained in (\ref{nonhomo}), the test presented in Theorem \ref{gof} for binary trees can be used for goodness-of-fit testing for conditioned Galton--Watson trees. The invariance property of the test statistic used in the binary tree setting is then interpreted as conditional on the choice of the number of leaves chosen to construct the $L$-tree: For a fixed number of leaves chosen to construct the $L$-tree, asymptotically, the invariance property in (\ref{invariance}) still holds. 
Starting with a tree $\tau_n$ with $n_l$ number of leaves, the requisite one-sample test invariant to the permutation of leaves, is constructed as follows: (1) construct a consistent estimator $\hat{\sigma}^2 $  of $\sigma^2$; (2) for each tree, construct an $L$-tree by randomly choosing a subset of the leaves; (3) using the test statistic defined as the product of $\hat{\sigma}^2 $ and the sum of the branch lengths of the $L$-trees, and Slutsky's theorem, construct an asymptotic rejection region based on the $\chi^2$ distribution. The invariance property is interpreted conditional on the number of leaves chosen randomly. The extension to the two-sample case is straightforward. 
\begin{theorem}\label{gof2}
	Suppose $\boldsymbol{\tau_n}=(\tau_{n_1},\ldots,\tau_{n_p})$ and $\boldsymbol{\eta_m}=(\eta_{m_1},\ldots,\eta_{m_q})$ are independent samples of conditioned Galton--Watson trees from $\pi_{\tau}$ and $\pi_{\eta}$ respectively with $\sigma_{\tau}^2$  and $\sigma_{\eta}^2$ as offspring variances, with respective consistent estimators  $\hat{\sigma_{\tau}}^2$ and  $\hat{\sigma_{\eta}}^2$. 
	%Suppose that $w_{i}$ and $u_{j}$ are the normalised distances from the root of randomly chosen vertices from the two samples and let $\hat{\sigma_{dp}}^2=2p(\sum_{i=1}^pw_i^2)^{-1}$ and $\hat{\sigma_{dq}}^2=2q(\sum_{i=1}^qu_i^2)^{-1}$. 
	Let $\mathbf{k}=(k_1,\ldots,k_p)$ be subsets of the leaves of $\tau_{n_i}$, where $k_i$ is chosen according to a uniform distribution on the leaves of $\tau_{n_i}$ and let $|k_i|$ denote the cardinality of set $k_i$.
	\begin{enumerate}
		\item
		For a fixed $\mathbf{k}$, define the critical function
		\begin{equation*}
		\phi(\mathbf{k},\alpha)=\left\{
		\begin{array}{l l}
		1 & \quad \text{if  }  \frac{1}{\hat{\sigma_{\tau}}^2} \sum_{i=1}^p s_i >\chi_{1-\alpha,2\sum_{i=1}^p|k_i|} ; \\
		0 & \quad \text{otherwise},
		\end{array}
		\right .
		\end{equation*}
		where $s_i$ are the total path lengths of $L(\tau_{n_i},k_i)$. Then, given $\mathbf{k}$, for the pair of hypotheses $H_0:\pi_{\tau}=\pi$ vs $H_1:\pi_{\tau}\neq \pi$, where $\pi$ is the density of a conditioned Galton--Watson tree,  the test given by $\phi(\mathbf{k},\alpha)$ is such that as $n_i \to \infty$, $E_{H_0}\phi(\mathbf{k},\alpha) \to \alpha$.
		\item 
		Choose  $\mathbf{g}=(g_1,\ldots,g_q)$ as the subset of leaves from $\boldsymbol{\eta_m}$ in a similar manner with cardinality $|g_j|$. 
		Then, the critical function
		\begin{equation*}
		\psi(\mathbf{k,g},\alpha)=\left\{
		\begin{array}{l l}
		1 & \quad \text{if  }  \frac{\hat{\sigma_{\eta}}^2\sum_{i=1}^p s_i }{\hat{\sigma_{\tau}}^2\sum_{j=1}^q r_j }>\left(\frac{\sum_{i=1}^p|k_i|}{\sum_{j=1}^q|g_j|}\right)F_{1-\alpha,2\sum_{i=1}^p|k_i|,2\sum_{j=1}^q|g_j|}; \\
		0 & \quad \text{otherwise} ,
		\end{array}
		\right .
		\end{equation*}
		where $r_j$ are the total path lengths of $L(\eta_{m_j},g_j)$,  defined, assuming no loss of generality that the numerator exceeds the denominator, for testing $H_0:\pi_{\tau}=\pi_{\eta}$ against $H_1:\pi_{\tau}\neq\pi_{\eta}$. The critical function $\psi$ is such that as $\min(n_i,m_j) \to \infty$ for every $i,j$, $E_{H_0}\psi(\mathbf{k,g},\alpha) \to \alpha$.
	\end{enumerate}
\end{theorem}

\subsection{Consistent estimation of offspring variance}
%\subsubsection{Estimation using Dyck paths}
For binary trees, under the Poisson process model, the parameter $\theta$ is estimated quite easily using the principle of maximum likelihood (using the density in (\ref{nonhomo})), since the size of the tree remains fixed. Under the conditioned Galton--Watson tree model, for a tree $\tau_n$, the density (\ref{LCAdensity}) is the asymptotic density of an $L$-tree with $k$ leaves, as $n \to \infty$. It is difficult to examine the behaviour of an estimator of $\sigma^2$, which is intricately dependent on $k$. It would be desirable to construct an estimator that does not depend on the number of leaves chosen to construct the $L$-tree. This can be achieved through the characterization of the CRT using a mapping of trees to a function space, known as \emph{Dyck paths}. 

\begin{figure}[!ht]
\centering
  \includegraphics[scale=0.6]{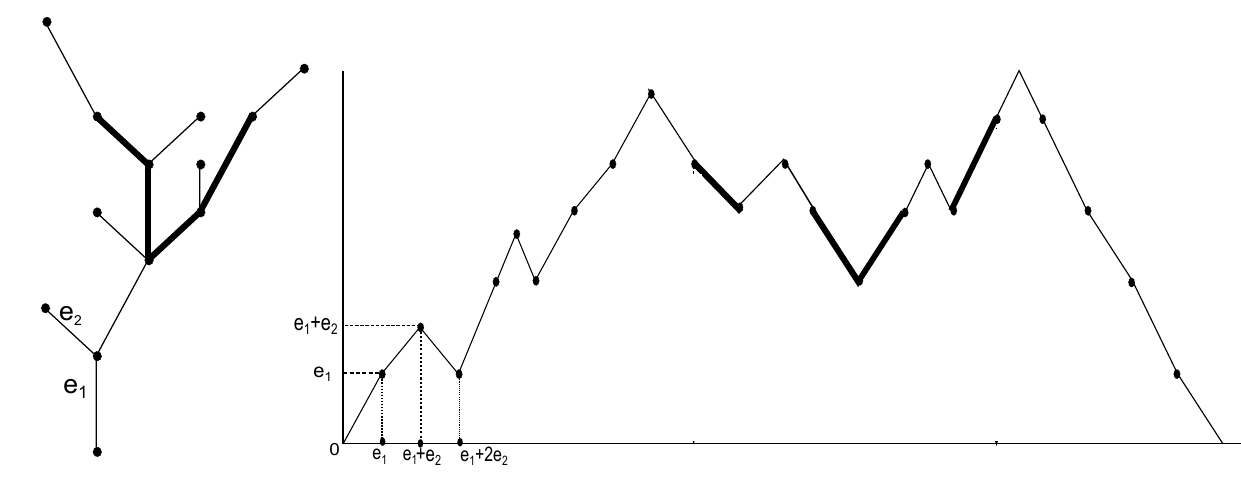}
\caption{\footnotesize A tree with root at the bottom and its corresponding Dyck path. The $x$ axis ranges from 0 to twice the sum of lengths of the edges; the Dyck path is constructed by traversing the tree in a depth-first manner at unit speed.}
  \label{figure1}
\end{figure}

Any rooted ordered tree of $n$ vertices can be uniquely coded by a traversal of the tree. When the traversal is a depth-first walk, one can construct a function, referred to as a Dyck path, that is bijective to the tree in the following manner. Imagine the motion of a particle that starts at time $t = 0$ from the root of the tree and then explores the tree from the left to the right, moving
continuously along the edges at \emph{unit speed} until all the edges have been explored and the particle has come back to the root. Each edge will be crossed twice in this evolution, hence the total time needed to explore the tree is $2l_n$, where $l_n$ is the total path length, or sum of the branch lengths. The walk can be represented as the value $H_n(s)$ of a continuous function $H_n: [0,2l_n] \rightarrow \mathbb{R}_{\geq 0}$ at time $s \in [0,2l_n]$ such that $H_n(s)=d(\text{root},v)$ where $v$ is the vertex obtained during the walk such that the sum of the edges traversed until $v$ is $s$, and $d(u,v)$ is the length of the unique path from vertex $u$ to vertex $v$. Figure \ref{figure1}, taken from \cite{pitman}, offers a more intuitive description. The map from $\tau_n$ to its Dyck path is a bijection. The Dyck path approach to analysis of tree-structured data was adopted by \cite{Shen-Marron}. Our interest in the Dyck path representation is captured in the following result for conditioned Galton--Watson trees on $n$ vertices; thus branch lengths are unit length with corresponding Dyck path $H_n:[0,2n]\to \mathbb{R}_{\geq 0}$. 
 \begin{theorem}\citep{aldous2}
\label{excursion}
Let $\tau_n$ be a conditioned Galton--Watson tree conditioned with an offspring distribution with mean $1$ and variance $\sigma^2 < \infty$. Let $H_n(k), 0 \leq k \leq 2n$ be the Dyck path associated with $\tau_n$. Then, as $n \to \infty$,
\[
\left\{\frac{1}{\sqrt{n}}H_n([2nt]), 0 \leq t \leq 1 \right\} \Rightarrow \left\{\frac{2}{\sigma}B^{ex}_t: 0 \leq t \leq 1\right\},
\]
where $B^{ex}$ is the standard Brownian excursion, and $\Rightarrow$ denotes weak convergence in $C[0,1]$.
%and $\Rightarrow$ implies weak convergence of processes in $C[0,1]$.
% the space of continuous functions on $[0,1]$, and $[\cdot]$ stands for the integer function.
\end{theorem}
Theorem \ref{excursion} implies an invariance principle with the Brownian excursion as the weak limit \emph{regardless} of the offspring distribution as long as it is critical with finite variance. Proposition \ref{binaryp} extends the result to certain sub-and super-critical cases. From this, asymptotic properties of the functionals of conditioned Galton--Watson tree models can be examined. The following Theorem presents a consistent estimator of $\sigma^2$. Recall that $d(u,v)$ represented the length of the unique path connecting vertices $u$ and $v$ in a tree. 
%For example, the distributions of the height of the tree and number of vertices at a given distance from the root can be approximated by the corresponding Brownian excursion functionals. One such functional, the distance from the root of a randomly chosen vertex, is used to construct an estimator of $\sigma^2$. 
\begin{theorem}\label{coro2}
Let $\boldsymbol{\tau_n}=(\tau_{n_1},\ldots,\tau_{n_p})$ be a random sample of conditioned Galton--Watson trees. 
% Consider a random sample of conditioned Galton--Watson trees $\boldsymbol{\tau_n}=(\tau_{n_1},\ldots,\tau_{n_p})$  from $\pi_{\sigma^2}$.
\begin{enumerate}
\item On each $\tau_{n_i}$, suppose $v_i$ is a vertex chosen according to a uniform distribution on $\mathcal{V}(\tau_{n_i})$. Then, the random variable $n_i^{-1/2}d(\text{root},v_i) \overset{d} \to W$, as $n_i \to \infty $, where $W$ is a Rayleigh distributed random variable with scale $1/\sigma$. 
\item Let $W_i$ be random variables that denote the normalised distance of a randomly chosen vertex from the root. Then, as $n_i \to \infty$ for each $i$, then $\hat{\sigma_d}^2=2p\left(\sum_{i=1}^pW_i^2\right)^{-1}$ is a consistent estimator of $\sigma^2$.
\end{enumerate}
\end{theorem}
%\noindent To summarize, the goodness-of-fit tests for conditioned Galton--Watson trees are constructed in the following fashion.
%\begin{enumerate}
%	\item Randomly choose a vertex in each tree, note its distance from the root, and compute an estimator of $\sigma^2$ using Theorem \ref{coro2}.
%	\item For each tree $\tau_n$, randomly choose a subset $B$ of its leaves, and construct the $L$- tree $L(\tau_n,B)$.
%	\item Compute the sum of the total path lengths of each constructed  $L$- tree, and compute the test statistic and construct the rejection region of the test using the result from Theorem \ref{gof2}.
%\end{enumerate}
%%%%%%%%%%%%%%%%%%

%%%%%%%%%%%%%%%%%%%%%%%
\section{Numerical illustrations}
\label{simulations}
\subsection{Simulations}
We use an efficient method to simulate conditioned Galton--Watson trees by employing the algorithm provided by \cite{LD} with a linear expected time. This enables us to simulate a large number of conditioned Galton--Watson trees, each of which contains a large number of vertices, and each tree is generated in expected linear time. Details of the simulation procedures can be found in Section 3 of Supplementary Material. Table \ref{one-sample-gof} reports the performances of tests in Theorem \ref{gof2} and a competing permutation test with the same test statistic. Rejection rates were computed by averaging over multiple permutations of the chosen leaves. Geo(0.5) denotes a Geometric distribution with probability 0.5; Bin(2,0.35) denotes subcritical conditioned Galton--Watson tree from a Binomial distribution with 2 trials and success probability 0.35; GW-Bin(2,0.5) denotes unconditioned Galton--Watson trees with Bin(2,0.5) distribution; Phylo.bd and Phlyo.coal correspond to phylogenetic trees based on birth-death processes on fixed taxa with speciation rate 2 and Kingman's coalescent process (see Section 1 of Supplementary Material for definition), respectively. The poor power against trees generated from a coalescent process is due to their connection with the CRT \citep{CH}, and this will be exploited in the data application. In Section 4 of Supplementary Material, we provide detailed results from examinations of the asymptotic behavior of test statistics, estimator of the offspring variance, and performances of the tests on binary trees from Poisson model. The results largely corroborate the theoretical findings.  

\begin{table}
	\small
	\begin{tabular}{|c|cc|cc|cc||cc|cc|cc|}
		%\hline
		\multicolumn{7}{c}{One-sample}&\multicolumn{6}{c}{Two-sample}\\
		\hline
		\multirow{2}{*}{Distribution}& \multicolumn{2}{c|}{$N=10$}&\multicolumn{2}{c|}{$N=100$}&\multicolumn{2}{c||}{$N=1000$}& 
		\multicolumn{2}{c|}{$N=10$}&\multicolumn{2}{c|}{$N=100$}&\multicolumn{2}{c|}{$N=1000$}\\
		%\hline
		&$\chi^2$& perm&$\chi^2$& perm&$\chi^2$& perm&$F$& perm&$F$& perm&$F$& perm\\
		\hline
		Geo(0.5)&0.09&0.15&0.05&0.08&0.03&0.09& 0.11&0.10&0.13&0.08&0.04&0.02\\
		Bin(2,0.5)&0.13&0.08&0.04&0.03&0.01&0.01& 0.08&0.14&0.08&0.09&0.03&0.06\\
		Bin(2,0.35)&0.10&0.16&0.12&0.07&0.06&0.08&0.21&0.14&0.13&0.10&0.04&0.01\\
		GW-Bin(2,0.5)&0.78&0.91&0.91&0.97&0.99&1.00&0.82&0.87&0.88&0.93&0.95&1.00\\
		Phylo.bd&0.81&0.83&0.89&0.91&0.98&0.94&0.92&0.83&0.96&0.97&0.99&1.00\\
		Phylo.coal&0.26&0.37&0.14&0.21&0.11&0.08&0.37&0.28&0.23&0.18&0.19&0.17\\
		\hline
	\end{tabular}
	\vspace*{4mm}
	%\end{table}
	%\vspace*{4mm}
	\caption{\footnotesize{One-sample: Rejection rate of  one-sample $\chi^2$ test, at $\alpha=0.01$, in Theorem \ref{gof2} and  permutation test with conditioned Galton--Watson trees containing 1000 vertices under sample sizes $n$, using $L$-trees constructed from 25 randomly chosen leaves. Two-sample: Rejection rate of two-sample goodness-of-fit F and permutation tests, having conditioned Galton--Watson trees with offspring distribution Bin(2,0.5) under $H_0$ against alternatives in column 1. }}
	\label{one-sample-gof}
\end{table}
\subsection{Data application: Detection of tumor heterogeneity using magnetic resonance images}
\label{data_application}
We illustrate the utility of the proposed tests through a novel approach to detecting tumour heterogeneity in brain cancer by constructing binary trees obtained from Magnetic Resonance (MR) images. 
%We first detail the dataset, the pre-processing steps and the objective of the study; we then briefly describe the relationship between the CRT and the trees obtained from hierarchical clustering, and present results from the tests. 

\begin{figure}[!ht]
\centering
\includegraphics[scale=0.7]{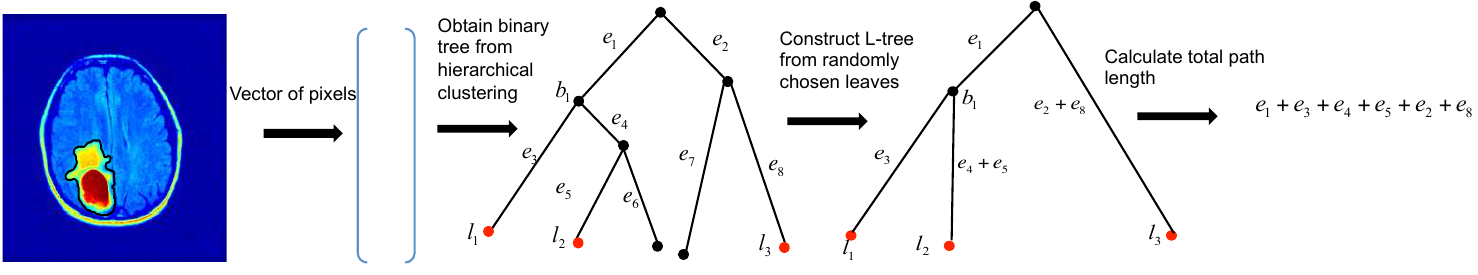}
\caption{\footnotesize Preprocessing steps starting from a MR image to the computation of total path length of an $L$-tree constructed from the resulting binary tree.}
\label{preprocessing}
\end{figure}

\noindent \underline{\emph{Data structure, pre-processing and key scientific question}}.  In this study, we used presurgical, T1-weighted post-contrast and T2-weighted/FLAIR images of 82 patients (26 women and 56 men) with histologically confirmed glioblastoma multiforme (GBM)---an aggressive form of brain cancer---from The Cancer Genome Atlas (TCGA) database. The images were downloaded from The Cancer Imaging Archive  at \texttt{https://www.cancerimagingarchive.net/} and are publicly available. 
We pre-processed the MR images and obtained 3-dimensional (3D) tumour volumes. Specifically, the images were registered spatially, followed by intensity bias correction using Medical Image Processing Analysis and Visualization software (MIPAV v6.0.0). The tumor region was segmented semi-automatically in 3D using the Medical Image Interaction Toolkit (\texttt{MITK.org}). Tumor regions were defined as a combination of the T1-contrast enhancing region as well as the FLAIR  hyperintense region, or specifically, the regions common to the T1-enhancing signal and FLAIR hyperintensity; this captures a combination of the tumor's enhancing component as well as the infiltrative edema component. The T1-weighted post-contrast and FLAIR tumor regions were delineated separately by a qualified neurosurgeon using the Medical Image Interaction Toolkit. The in-plane resolution of the image was $1 \text{mm} \times 1 \text{mm}$. Our analysis was based on T2-weighted intensities from \emph{only} the segmented regions; see image in Figure \ref{preprocessing} where the segmented region is outlined in black. 

Tumor heterogeneity expressed though pixel intensities indicates a latent ordering of groups of pixels with similar intensities that represent similar etiologies. The problem of interest is to appropriately characterize tumor heterogeneity in the brain tumors. Specifically,  we consider detection of tumor heterogeneity between patients with GBM who have long survival times versus those who have short survival times. 

\noindent \underline{\emph{Methods}}. Current approaches to this problem are based on simple summaries of the entire image such as skewness or kurtosis of the probability density of the intensities, which fail to take into account the structural complexities of the pixel-level intensities, while also neglecting spatial information; see \cite{NJ} for a detailed review. In contrast, our approach is based on exploring the clustering properties inherent in the density, extracted through the number and size of modes. Our approach, illustrated in Figure \ref{preprocessing}, is as follows.\\
1. For each patient, from a single axial slice, obtain the pixels from the segmented tumour.\\
2. Construct a binary tree/dendrogram by implementing an agglomerative hierarchical clustering algorithm on the pixels.\\
3. Randomly choose a subset of the leaves, construct the corresponding $L$-tree, and compute its total path length.

On the tree, the individual intensities appear as leaves, clusters of similarly valued intensities are the internal vertices, and branch lengths represent distances between the clusters. Such a representation expresses heterogeneity in the tumor image though the clustering of intensity values. A rationale for using vectorized pixel intensities, disregarding the spatial information, is that when the chief objective is to characterize or classify tumor images based on heterogeneity, the location of the heterogeneous regions on the image is of little relevance: the classifying procedure needs to only \emph{detect} the heterogeneity in the entire image in relation to another image. Our view of heterogeneity is closely linked to the density-cluster tree \citep{hartigan, AT} from density-based clustering, wherein the focus is on the \emph{multiple modes within the probability density} of the intensity values: each branch represents the high-density clusters within a single mode of the density, i.e. groupings of pixels that are not separated easily. Branch points of the tree represent values at which a new mode of the density, or a new cluster, emerges. Hence, on the scale of the data at multiple resolutions, the branch lengths are indicative of how long a particular mode within the density lasts before being broken up; the sum of the branch lengths is the corresponding cumulative measure on the scale of the data, capturing heterogeneity through the number and size of modes in the density of the pixel intensities.

\begin{figure}[!htb]
\begin{tabular}{cc}
\includegraphics[trim=105mm 45mm 0mm 60mm,clip,scale=0.4]{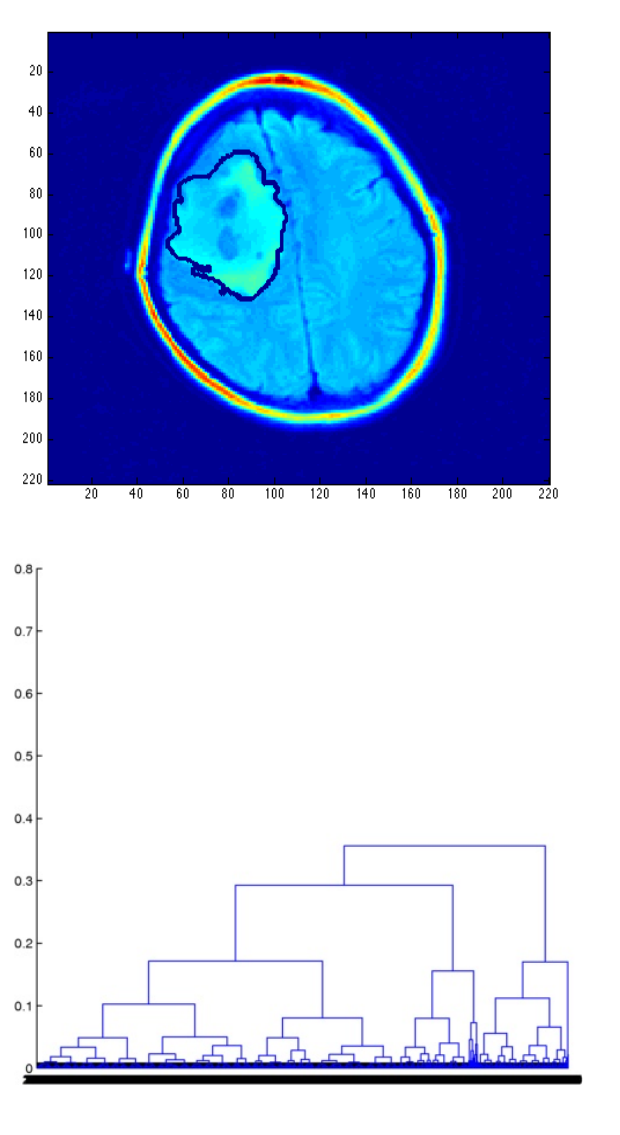}&
\hfill
\includegraphics[trim=105mm 45mm 0mm 60mm,clip,scale=0.4]{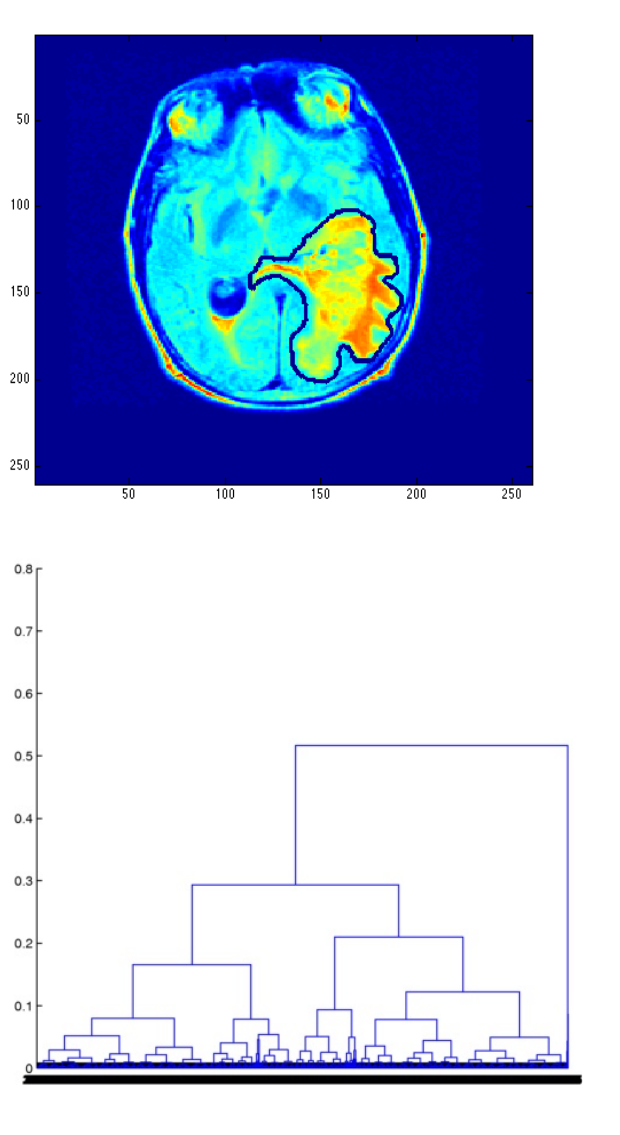}
\end{tabular}
\caption{\footnotesize T2-weighted Brain MR images of segmented tumors (outlined in black) and dendrograms for a patient who experiences (left) a long survival time of 57.8 months and (right) a short survival time of 0.723 months. Dendrograms were constructed from pixels obtained from only the segmented regions.}
\label{tumours}
\end{figure}

Subsequently, the test in Theorem \ref{gof2}, with the test statistic that relates the pixel-cluster distances to the branch lengths of an $L$-tree, can be used to detect group differences between patients with GBM who experience long survival times ($>$ 12 months) and those who experience short survival times ($\leq$ 12 months). As an illustration, consider the images of the patients who respectively correspond to the short (0.723 months) and long (57.8 months) survival times, as shown in the top panel of Figure \ref{tumours}. Their respective trees  have total path length 3.69 and 2.13 respectively. Intuitively, the MR image for the patient with the shorter survival time, with the tumor appearing to be in an advanced stage in the image on the right in Figure \ref{tumours}, should have pixel intensities with richer and varied clustering tendencies, which are evident in the form of more branches with smaller branch lengths. 
From the analysis of the rejection rates of the tests in Table \ref{one-sample-gof}, we note that the two proposed test techniques do not have good power against ultrametric trees such as coalescent trees. These trees are characterised as possessing leaves that are equidistant from the root. Dendrograms arising from agglomerative hierarchical clustering are ultrametric, and a surprising connection exists, in the asymptotic regime, between the dendrograms and the $L$-trees that characterize the CRT;  see \citet{CH} for details. 

\subsection{Results}
We use agglomerative hierarchical clustering in \texttt{MATLAB} with Single linkage function to obtain dendrograms of the image intensities. \cite{CM} noted that the single linkage function is the most stable choice for hierarchical clustering in the sense that the distance is small between the original tree and the tree obtained after a small perturbation of the input data. Similar constancy results under Single linkage function are available for the density-cluster \citep{hartigan}. However, the large sizes of the images considered in our setting led to trees that did not appear significantly different for different choices of the linkage functions; this was corroborated by the results of the $F$ and permutation tests, which were unchanged in relation to different choices of linkage function. In order to use the test prescribed in Theorem \ref{gof2}, we estimate the variance $\sigma^2$ with $\hat{\sigma_d}^2$ from Theorem 5. The small sample sizes---31 and 42 patients in groups who experienced short and long survival times, respectively---do not guarantee a consistent estimate of $\sigma^2$; however, this does not appear to affect the efficacy of the tests. 

\begin{table}
\small
\begin{tabular}{|cccc|ccc|}
\hline
\multicolumn{4}{|c|}{\it Long vs short survival }&\multicolumn{3}{|c|}{\it Short vs short survival}\\
\cline{1-7}
 $\%$ of leaves  & Null dist. & F-statistic & Decision at 1\%& Null dist & F-statistic & Decision at 1\% \\
\hline
10&  $F_{\scriptscriptstyle{11638,11734}}$&0.780&Reject&$F_{\scriptscriptstyle 91494,118142}$&0.990&Do not reject\\
20&  $F_{\scriptscriptstyle 23310,23515}$&0.775&Reject&$F_{\scriptscriptstyle 61721, 74616}$&0.985&Do not reject\\
30&  $F_{\scriptscriptstyle34793,35286}$&0.774&Reject& $F_{\scriptscriptstyle 66558, 91285}$&0.997&Do not reject\\
40&  $F_{\scriptscriptstyle46648,47070}$&0.772&Reject& $F_{\scriptscriptstyle 88392,101228}$&1.001&Do not reject\\
%50&  $F_{\scriptscriptstyle291375,293975}$&0.776&Reject\\
\hline
\end{tabular}
\vspace*{3mm}
\caption{\footnotesize  At different values for the percentage of leaves chosen to construct the $L$-trees of the dendrograms from the images, the results of the $F$ test from Theorem \ref{gof2} to test the null hypothesis of equality of distribution of the two groups. }
\label{test-results}
\end{table}

Using the total path length, scaled by the estimate of the variance, of an $L$-tree constructed by choosing various percentages of leaves at random as the test statistic, the $F$ test from Theorem \ref{gof2} rejected the null hypothesis of equality of distribution, at 1\% significance level, between the groups with long versus short survival times, and failed to reject $H_0$ when the labels were randomly permuted. Within each group, when two subgroups were randomly chosen and compared for group differences, the test failed reject the null hypothesis of equality of distribution at 1 \% significance level. Details are reported in Table \ref{test-results} for the comparison of groups with long versus short survival times, and within groups with survival times. Note that the choice of the proportion of leaves used to construct the $L$-trees does not seem to affect the conclusions of the test. This is particularly encouraging since computation of $L$-trees from observed trees can be computationally expensive. The permutation test provided identical results for the same test statistic at 2000, 5000 and 10000 permutations. Importantly, the conclusions of the tests were unaltered when the linkage function was changed. Scientifically, this points to the existence of significant GBM tumor heterogeneity among patients with different durations of survival.
%%%%%%%%%%%%%
\section{Discussion}
\label{discussion}
%Tree-structured data can be regarded as data observed on acyclic graphs with a distinguished vertex that is referred to as the root. This is in contrast to networks wherein cycles are allowed and are of chief interest. This precludes the use of a large number of statistical methods currently available for network data, when dealing with tree-structured data.
Aldous' papers on the CRT and variants (see, in addition, \cite{aldous4,aldous5}) provide useful distributional results and connections to common stochastic processes, which in principle can be harnessed in modeling tree functionals (for e.g. total height, Wiener index) by corresponding functionals of the limit stochastic processes.  
%One immediate extension to this work is to approximate the distributions of local tree functionals by the corresponding Brownian excursion functionals. For example, the Wiener index of a tree, popular in phylogenetics and chemistry, is exactly $\frac{2}{n}A_n$ where $A_n$ is the area under the curve of the Dyck path of a tree with $n$ vertices. The distribution of $A_n$ can be approximated by the distribution of the Brownian excursion area, which is well known (albeit difficult to compute); see \cite{janson} and references therein for details.
However, the absence of information on the rates of convergence of conditioned Galton--Watson trees to the CRT prevents us from obtaining a clear idea as to the size of the trees that assure the validity of our results. 

Apart from the need to retain the projective property, another reason for using the multiplicative model $\lambda(t)=\theta t$ for the intensity function is to model the dependence of the CRT on the variance $\sigma^2$ of the offspring distribution of the conditioned Galton--Watson trees. The parameter $\sigma^2$, as can be seen from the construction and Lemma 1, appears naturally in the limiting density of the $L$-trees as a scale parameter; this is very useful for its consistent estimation using a maximum liliehood estimate, described in Theorem 5. 

An approach for testing not explored in this article is based on using the normalized Dyck paths to define likelihoods. Recall that $\frac{2}{\sigma}B^{ex}$ is the limit of normalized Dyck paths that uniquely code conditioned Galton--Watson trees.  The CRT is connected to $\frac{2}{\sigma}B^{ex}$ in the following manner: Consider uniform order statistics $U_{1:k}< \cdots<U_{k:k}$ and set $V_{i}=\min_{U_{i:k}\leq t \leq U_{i+1:k}}\frac{2}{\sigma}B^{ex}(t)$. 
%Draw an edge of length $\frac{2}{\sigma}B^{ex}(U_{1:k})$ and label one end as the root and the other end as $U_1$. Inductively, from $U_{i:k}$ move back a distance
%$\frac{2}{\sigma}B^{ex}(U_{i:k})-V_i$ toward the root, draw a new edge of length $\frac{2}{\sigma}B^{ex}(U_{i+1:k})-V_i$ and label the new endpoint $U_{i+1}$.  The resulting binary tree on $k$ vertices with $k-1$ edges has the density given in (\ref{LCAdensity}). For a tree $\tau_n$ with $n$ vertices, let $0 =U_{0:n}< U_{1:n}<\cdots<U_{n+1:n}=1$ be uniform order statistics and let  $V_{i}=\min_{U_{i:n}\leq t \leq U_{i+1:n}}\frac{2}{\sigma}B^{ex}(t)$. 
Then the $(2n+2)$-dimensional vector taking values in $\mathbb{R}_+^{2n+2}$ as
$X_n=\Big(\frac{2}{\sigma}B^{ex}(U_{i:n}),\frac{2}{\sigma}B^{ex}(V_i)\Big)$ has the same distribution as
%\end{equation*}
%Based on this construction, it can be seen that the distribution of the random vector $X_n$ defines a distribution on the random tree constructed with $n$ vertices. One way of looking at the density in (\ref{LCAdensity}) via this construction is as the density of the random variable that is the total variation of the function obtained via a linear interpolation between points in $X_n$. Using this approach it was shown in \cite{pitman2} that
$\frac{\sigma(2\Gamma_{n+1})^{1/2}}{4}\Big(U_{i-1:n}-V_{i-1},U_{i:n}-V_{i-1}; 1 \leq i \leq n+2\big|\cap_{i=1}^n(U_{i:n}>V_i) \Big),$
where $U_{n+2:n}:=1$ and $\Gamma_{n+1}$ is a Gamma random variable with shape $n+1$ and scale 1\citep{pitman2}. Although in principle, it would be reasonable to define a parametric class, the distribution of $X_n$ is not easy to compute.

The data application considered in this article represents, to our knowledge, the first attempt at characterizing tumor heterogeneity from images that use tree representations. However, important extensions such as establishing correspondence between images, incorporating spatial information of the pixels, using covariate information for each patient, and developing methods for images obtained from a  longitudinal study of the patients, are part of our current work. Much remains to be done in this direction.
\section{Acknowledgements}
We thank Joonsang Lee, Juan Martinez, Shivali Narang and Ganesh Rao for their assistance with MR images in the data application, and the referees for helpful comments. 
%%%%%%%%%%%%%%%%%%%%%
\bibliography{ref}
\bibliographystyle{plainnat}
\newpage
%%%%%%%%%%%%%%%%%
\begin{center}\textbf{SUPPLEMENTARY MATERIAL}\end{center}
\setcounter{section}{0}
\section{Definitions}
%\subsection{Definitions}
\noindent \underline {\emph{Coalescent process}}:
Introduced by \cite{kingman}, the coalescent process can be described as follows. Fix $k$; at time $t=0$ there are $k$ individuals $(1,\ldots,k)$; at each time $t>0$, there is a partition of the individuals in to $j$ clusters , $1 \leq j \leq k$. The process evolves according to the rule: in time $[t,t+dt]$, each of the ${j \choose 2}$ pairs of clusters with probability $dt$ coalesce (enjoin) into a single cluster. A tree constructed using such a process is not a conditioned Galton--Watson tree. \\
\underline{\emph {Projective probability measures}}. 
Let $T$ be a finite or countable index set. For any collection of spaces $S_t$ with $t \in T$, for $I \subset T$ let $S_I=\times_{t \in I} S_t$. A family of probability measures $\mu_I, I \in T$ is projective (p. 92 of \cite{Kallenberg}) if
 $$\mu_J(\cdot \times S_{J\setminus I})=\mu_I, \quad I \subset J \text{ in } T .$$
\underline{\emph{Distinguishable hypotheses}}. 
 Suppose $\Theta$ is an index set and $\Theta_0$ and $\Theta_1$ are disjoint subsets of $\Theta$ such that $\Theta_0 \cup \Theta_1=\Theta$. Denote by $H_0$ and $H_1$ the null and alternative hypotheses that $\theta$ is a member of either $\Theta_0$ or $\Theta_1$. Then, the set of probability measures $\Big\{P_{\theta}: \theta \in \Theta\Big\}$ is distinguishable (\cite{Rao}) if
 \begin{enumerate}[(i)]
 	\item $P_\theta\neq P_{\theta'}$ for all distinct $\theta,\theta' \in \Theta$; and
 	\item there is at least one Borel set $A$ such that $ P_\theta(A)\neq P_{\theta'}(A)$ for $\theta \in H_0$ and $\theta' \in H_1$.
 \end{enumerate}

\section{Technical details and proofs}
\subsection{Properties of $f_{k,\theta}$ on binary trees obtained from Poisson process model} If the family of probability measures that correspond to the densities $f_{k,\theta}$ is projective,  the binary tree with $k-1$ leaves, $\tau(k-1)$, obtained upon deleting an edge from $\tau(k)$ possesses the density $f_{k-1,\theta}$. This guarantees that when $f_{k,\theta}$ is parameterized by $\theta$, the interpretability of $\theta$ is retained as leaves are deleted or added, facilitating the development of consistent statistical procedures. Furthermore, in the interest of constructing tests, we require that the probability measures are also distinguishable. 
 \begin{theorem}\label{projective}
 	Suppose  $\mu_k$ is a family of probability measures on $\mathcal{T}_{2k} \times \mathbb{R}_{+}^{2k-1}$ corresponding to the densities $f_k$. 
 \begin{enumerate}[(a)]
 \item 
 The family $\mu_k$ is projective. The same holds for the sequence $\mu_{k,\theta}$, which corresponds to the sequence  $f_{k,\theta}$ for $\theta \in \Theta$.
\item 	The parametric family $\{\mu_{k,\theta}: \theta \in \Theta \subset \mathbb{R}^+\}$ is distinguishable.
\end{enumerate}
 \end{theorem} 
\noindent \textbf{Proof of} (a). The proof shows that the densities satisfy the projective property, which implies that the corresponding probability measures also satisfy this property (see p. 93 of \cite{Kallenberg}). The proof of the result for the parameterized case is applicable to the nonparameterized case with $\theta=1$. Hence, we assume that 1 is an element of $\Theta$. The density $f$ for a fixed $\theta$ is exchangeable with respect to the edge lengths since $s=\sum_{i=1}^{2k-1}x_i$ is invariant to permutations of $x_i$. This implies that relabeling the vertices and the edges leaves the distribution  unchanged. For ease of notation, let
 	$$\mathcal{L}_k=\mathcal{T}_{2k} \times \mathbb{R}_+^{2k-1} .$$
 	Consider $\tau(j) \in \mathcal{L}_j$ with $j$ leaves. The probability kernel $p_j$ from $\mathcal{L}_1 \times \cdots \times \mathcal{L}_{j-1}$ to $\mathcal{L}_j$
 	is defined in terms of the conditional density
 	$$f_{j,\theta}(\tau(j)  | \tau(j-1))\propto \frac{s'}{s}e^{-\frac{(s'^2-s^2)\theta}{2}},$$
 	with $s=x_1+\cdots+x_{2j-3}$ and $s'=s+x_{2j-2}+x_{2j-1}$. By induction on $j$, we can extend the existence of the probability kernel $p_k$ to $\mathcal{L}_k$ with conditional density
 	$$f_{k,\theta}(\tau(k)|\tau(k-1))= \frac{s'}{s}e^{-\frac{(s'^2-s^2)\theta}{2}},$$
 	where $s=x_1+\cdots+x_{2k-3}$ and $s'=s+x_{2k-2}+x_{2k-1}$. By Theorem 5.17 in \cite{Kallenberg}, we can assert the existence of the tree $\tau(k)$ with distribution $p_1 \otimes \cdots \otimes p_k$; in other words, the distribution on $\tau(k)$ can be defined via the conditional densities as
 	\begin{align*}
 	f_{k,\theta}(\tau(k))&=f_{1,\theta}(\tau(1))f_{1,2,\theta}(\tau(2)|\tau(1)) f_{2,3,\theta}(\tau(3)|\tau(2))
 	\ldots f_{k-1,k,\theta}(\tau(k)|\tau(k-1)).
 	\end{align*}
 	Straightforward computation with the conditional densities verifies this fact. 
\qed
\\
\textbf{Proof of} (b). 
Suppose $B$ is a Borel subset of $\mathcal{C}_k=\mathcal{T}_{2k} \otimes \mathbb{R}_+^{2k-1}$. Suppose we define a relation $\sim$ on subsets $B_1$ and $B_2$ of $\mathcal{C}_n$ as $B_1 \sim B_2$ if they contain \emph{all} trees with $2k$ vertices; by this we mean that the ``shape" of the tree is disregarded and imply that all trees with $2k$ vertices are equivalent. Note that $\sim$ is an equivalence relation and generates the quotient class $\mathcal{C}_{2k}^\sim=(\tau(k),(e_1,\ldots,e_{2k-1}))$ with $(e_1,\ldots,e_{2k-1}) \in \mathbb{R}_+^{2k-1}$ and $\tau(k)$ is the canonical tree with $2k$ vertices and $k$ leaves. The Borel sets of $C_{k}^\sim$ are the usual open rectangles generating the Euclidean space $\mathbb{R}_+^{2k-1}$. Note that the law $\mu_{k,\theta}$ assigns different mass to distinct elements in $\mathcal{C}_k^\sim$. We are hence interested primarily in Borel subsets of $C_k^\sim$ and restrict our examination of distinguishability to this equivalence class.

Suppose the null hypothesis $H_0$ is that $\theta \in \Theta_0$ and the alternative hypothesis $H_1$ is that $\theta \in \Theta_1$, where $ \Theta_0 \cup  \Theta_1= \Theta$ and $ \Theta_0\cap  \Theta_1=\emptyset$. Furthermore, the distribution function associated with $\mu_{k,\theta}$, for a tree $\tau_n$, corresponding to the density $f_{k,\theta}$ (continuous as a mapping $\theta \mapsto f_{k,\theta}$),
\[
F_{\theta}(x_1,\ldots,x_{2k-1})=\displaystyle \int_{0}^{x_1}\ldots\int_{0}^{x_{2k-1}} f_{\theta}(\tau(k))\quad de_{1}\ldots de_{2k-1}
\]
is continuous for each vector $(x_1,\ldots,x_n)$ representing edge lengths. Consider the function
$$g_k(x_1,\ldots,x_{2k-1})= \left[\prod_{i=1}^{k-1}\frac{1}{2i-1}\right]^{-1} (x_1+\cdots+x_{2k-1}).$$
Note that 
$$\displaystyle \int_{a_1}^{b_1}\ldots\int_{a_{2k-1}}^{b_{2k-1}} g_k(x_1,\ldots,x_{2k-1})\quad dx_{1}\ldots dx_{2k-1} < \infty$$
on any bounded rectangle $\prod_i[a_i,b_i]$, and hence is locally integrable. Since $f_{k,\theta} \leq g_k$ for every $\theta \in \Theta$, we can claim from Theorem 1 on p.7 of \cite{Rao}, that the sequence $\mu_{k,\theta}$ is distinguishable.
\qed\\
\begin{proposition}\label{gamma}
Suppose $\tau(k)$ is a binary tree with $k$ leaves and branch lengths $x_1,\ldots,x_{2k-1}$ generated from the non-homogeneous Poisson model. Then the total path length $\sum_{i=1}^{2k-1}x_i$ follows a Gamma distribution with shape $k$ and scale $2$.  
\end{proposition}
\begin{proof}
From the construction described it is easy to see that the distribution of the total path of $\tau(k)$ coincides with the distribution of the time of the $k$th event, say $S_k$. The density of $S_k$ for a  non-homogeneous Poisson process is 
$$f_{S_k}(t) = e^{-m(t)}m(t)^{k-1}\lambda(t)\frac{1}{(k-1)!} \quad t \geq 0,$$
where $m(t)$ is the intensity function. Noting that when $\lambda(t)=t$, $m(t)=t^2/2$, we obtain
$$f_{S_k}(t) = e^{-t^2/2}\left(\frac{t^2}{2}\right)^{k-1}\frac{t}{(k-1)!} \quad .  $$
Making a change of variable $y=t^2$, we see that $f_{S_k}(y)$ coincides with a Gamma density with shape $k$ and scale $2$. 
\end{proof}
\noindent \underline{\emph{Invariance to permutation of leaves}}. The permutational symmetry focuses attention on the class of test procedures that are invariant to the action of the finite symmetry group on the tree branch lengths. We formalize this notion in the following manner. Conditional on $k_i$, for a sample of trees $\tau(k_i), i=1,\ldots,N$, the joint data space is $\otimes_i^N\mathcal{X}_{k_i}$, for uniformly chosen topologies. Thus, each tree $\tau(k_i)$ is associated with its symmetric group of all permutations of its $2k_i-1$ branch lengths, where a permutation of length $2k_i-1$ is a bijection $\sigma:[2k_i-1]\to [2k_i-1]$, where $[2k_i-1]=\{1,2,\ldots,2k_i-1\}$. In order to define an invariant test procedure across an \emph{observed} sample of trees, we can define the group $G_{\text max}$ as the symmetric group of all possible permutations of the set $[2k_{\text{max}}-1]$ where $k_{\text{max}}:=\text{max}_i k_i$. For each tree $\tau(k_i)$ with $k_i \neq k_{\text{max}}$, we associate the set 
$$G_i:=\Big\{\sigma(1),\sigma(2),\ldots,\sigma(2k_i-1),2k_i,2k_i+1,\ldots,2k_{\text{max}}-1 | \sigma:[2k_i-1] \to [2k_i-1]\Big\}.$$
The set $G_i$ contains all possible permutations of the first $2k_i-1$ elements of $G_\text{max}$ and leave the rest unchanged, which leads to the following proposition.
\begin{proposition}\label{perm_group}
For each $i$, the set $G_i$ is a subgroup of $G_{\text{max}}$.
\end{proposition}
\begin{proof}
$G_\text{max}$ is a finite group of $(2k_{\text{max}}-1)!$ elements, and the cardinality of $G_i$ is also finite since $k_i <k_{\text{max}}$ for all $i$. The group operation of $G_\text{max}$ is composition of permutations. If $g^1_i$ and $g^2_i$ are two elements of $G_i$, then the action of $g^1_i$ on $\{1,2,\ldots,2k_\text{max}-1\}$ is $\{\sigma(1),\ldots,\sigma(2k_i-1),2k_i,,\ldots,2k_{\text{max}}-1\}$, where $\sigma(j)\in \{1,2,\ldots,2k_i-1\}$. The composition $g^2_i(g^1_i)$ is essentially the composition of $\sigma$ with itself, the result of which is an element in $\{1,2,\ldots,2k_i-1\}$. Hence $G_i$ is closed under composition. Using Proposition 2.69, p.149 of \cite{JR}, which states the conditions necessary for a finite subset of a group to be a subgroup, we can claim that $G_i$ is a subgroup of $G_\text{max}$ for every $i$.
\end{proof}

\noindent Any test function $\phi_\theta(\tau(k_1),\ldots,\tau(k_N))$ for testing the hypothesis that the sample of trees is from the generating model prescribed by $f_\theta$ can be reduced to the function $\phi_\theta(\mathbf{x}_i, i=1,\ldots,N)$ where $\mathbf{x}_i=(x_{1i},\ldots,x_{2k_i-1})$. Since for every action $g$ on $\mathcal{X}_k$, the induced action on the parameter set $\Theta$ is the identity mapping, the test function $\phi_\theta$ should satisfy, every $\theta$, 
\begin{equation}
\label{invariance}
\phi_\theta(\mathbf{x}_i, i=1,\ldots,N)=\phi_\theta(g^j_i(\mathbf{x}_i), i=1,\ldots,N), \quad g^j_i \in G_i, j=1,\ldots,(2k_i-1)!\quad.
\end{equation}
Let $f_\theta$ for $\theta \in \Theta \subset \mathbb{R}^+$ be the density of a tree obtained from the non-homogeneous Poisson model with intensity function $\lambda(t)=\theta t$, regardless of the number of leaves in the tree. From result (b) of Theorem \ref{projective}, we are guaranteed distinguishability. The goal is to construct a test function $\phi$ that satisfies (\ref{invariance}), with a desired $\alpha$ level, distinguishing between distinct values of $\theta$. The total path length statistic, that completely determines the density $f$ once the number of leaves is observed, will be used.\\
\\
\textbf{Proof of Theorem 1 in paper}.
 It is straightforward to compute the maximum likelihood estimate $\hat{\theta}$ of $\theta$ under $f_{\theta}$. We note that $\hat{\theta}=p(\sum_{i=1}^ps_i)^{-1}.$ \\
 \emph{One sample}: The likelihood ratio test statistic 
 $$\Lambda \propto \left(\frac{\theta'\sum_{i=1}^ps_i}{p}\right)^{\sum_{i=1}^pn_i}\text{exp}\left[-\frac{p\theta'\sum_{i=1}^ps_i}{2p}+\frac{p}{2}\right] .$$
Consider the function $h(t)=t^\alpha e^{1-t}$ for $t>0$ and $\alpha>0$. Then, $g'(t)=g(t)(\alpha/t-1)$, which is positive if and only if $\alpha >t$. Setting $t=(\theta'\sum_{i=1}^ps_i)p^{-1}$, and $\alpha=(2\sum_{i=1}^pk_i)p^{-1}$, the test rejects $H_0$ only for large values of $\sum_{i=1}^ps_i)$. Using Proposition \ref{gamma} and the reproductive property of independent Gamma random variables, we obtain the distribution of the test statistic. The critical function $\phi$ preserves permutational symmetry within the branch lengths of each tree since the test statistic is a sum of the total path lengths, each of which is invariant to permutations of the branch lengths.\\
\qed
%%%%%%%%%%%%%%%%%%%%%%%%%%%%%%%%%%%%%%%%%%%
%%%%%%%%%%%%%%%%%%%%%%%%%%%%%%%%%%%%%%%%%%%
\subsection{Properties of Conditioned Galton--Watson tree models}
This section collects some technical results that allow us to defined coherent models using conditioned Galton--Watson trees and their connection to the CRT, and proofs of results in the paper. \\
\textbf{Proof of Proposition 1 in paper}. 
\cite{DPK} proved that under certain conditions, the study of non-critical Galton--Watson processes reduces to the study of critical ones. Specifically, if $\lambda>0$ is a parameter such that $N_\lambda=\sum_{i\geq0}\pi_i\lambda^i <\infty$, we set $\pi_i^\lambda=\pi_i\lambda^i/N_\lambda$. Then a conditioned Galton--Watson tree with offspring distribution $\pi$ has the same distribution as a conditioned Galton--Watson tree with offspring distribution $\pi^\lambda$. As a consequence, if we can find a value of $\lambda$ such that the conditioned Galton--Watson tree under $\pi^\lambda$ becomes critical, then its asymptotic behavior is similar to that in the critical case. We examine each case separately and observe that there exists a finite $\lambda$ such that $N_\lambda<\infty$. 
\begin{enumerate}[(i)]
\item Note that 
$$N_\lambda=\sum_{i\geq1}(1-p)^{i-1}p\lambda^i=\frac{p\lambda}{1-\lambda(1-p)} $$
if $(1-p)\lambda<1$ for all $p$. Hence by choosing $\lambda<1/(1-p)$ we obtain the result with the requirement that $\lambda<1$ (only then is $\lambda<1/(1-p)$ for all $0<p<1$). Here, $\lambda=1/(2(1-p))$ results in a unit mean under the new distribution. 
\item Under this setup 
$$N_\lambda=(1-p)^2+2p\lambda(1-p)+p^2\lambda^2 $$
which is finite for all $\lambda<\infty$. The unit mean is attained at $\lambda=(1-p)/p$.
\item Here, $ N_\lambda=(1-p)+p\lambda^2,$
which again is finite for all finite $\lambda$; the unit mean is attained at $\lambda=\sqrt{(1-p)/p}$.
\item For this case,
$N_\lambda=p_0+p_1\lambda +(1-p_0-p_1)\lambda^2$, which is finite for all $\lambda<\infty$ with unit mean being attained at $\lambda=p_0/(1-p_0-p_1)$.
\end{enumerate}
%For $\pi_0=1-p$ and $\pi_2=p$ with $0<p<1$, $N_\lambda=(1-p)+p\lambda^2 <\infty$ for all $\lambda <\infty$. This implies that 
%$$\pi_0^\lambda=\frac{(1-p)}{(1-p)+p\lambda^2}\quad \text{and}\quad \pi_2^\lambda=\frac{p\lambda^2}{(1-p)+p\lambda^2} .$$
%If $\xi_\lambda$ is a random variable from $\pi^\lambda$, then 
%$$E\xi_\lambda=\frac{2p\lambda^2}{(1-p)+p\lambda^2} .$$ Setting $E\xi_\lambda=1$ and solving for $\lambda$, we obtain the result. 
\qed
\\
%%%%%%%%%%%%%%%%%%%%%%%%%%%%%%%%%%%%%%
\underline{\emph{Interpretation of $\sigma^2$ in $L$-tree model}}. While the offspring variance $\sigma^2$ has a clear enough interpretation in the context of conditioned Galton--Watson processes, its role in the density of $L$-trees is not clear. The specific question is the validity of $\{f_{k,\sigma^2}\}$ as a parametric class for $L$-trees of conditioned Galton--Watson trees with possible extension to the conditioned Galton--Watson trees themselves: The parameterized distribution on the $L$-tree needs to  be extended to $\tau_n$ for every $n$ in order to retain the interpretability of $\sigma^2$. The issue of \emph{extendability} is crucial. Suppose $\tau_n$ on $n$ vertices contains $n_l$ leaves and an $L$-tree with $k<n_l$ leaves is chosen randomly as described in Section 2 of the main document. With leaves as $l_i$, clearly $L(\tau_n, \{l_1,\ldots,l_{n_l}\} )=\tau_n$ and it is necessary that any $f_{k,\sigma^2}(\cdot)$ defined on $L(\tau_n,B)$ for any $B \subset \{l_1,\ldots,l_{n_l} \}$ with $|B|=k$ can be extended to $\tau_n$ upon inserting leaves from $\{l_1,\ldots,l_{n_l }\}-B$ to $L(\tau_n,B)$, while retaining the interpretability of $\sigma^2$.  
 
Recall that $\tau_n=(\mathcal{V}(\tau_n),\mathcal{E}(\tau_n))$ resides in $\mathcal{T}_n \times \mathbb{R}_+^{n-1}$; its least common ancestor tree corresponding to $B \subset \{l_1,\ldots,l_{n_l}\}$, lies in an appropriate subspace $\mathcal{T}_{2k} \times \mathbb{R}_+^{2k-1}$ with $2k<n$. Then, for all $k$ and  $n$ with $2k<n$, the parametric density defined on $\mathcal{T}_{2k} \times \mathbb{R}_+^{2k-1}$ can be extended to $\mathcal{T}_n \times \mathbb{R}_+^{n-1}$, or is \emph{$n$-extendable}, if $f_{k,\sigma^2}(\cdot)$ on $L(\tau_n,B)$ can be recovered by marginalization over $f_{\sigma^2}(\tau_n)$ on $\mathcal{T}_n \times \mathbb{R}_+^{n-1}$ for every $\sigma^2$.  The following result states that the density in (2) of the main document can be used as a consistent statistical model satisfying the extendability of $\sigma^2$.
\begin{proposition}
	The class $\Big\{f_{k,\sigma^2}: \sigma^2 \in \mathcal{S}\Big\}$ on $\mathcal{T}_{2k} \times \mathbb{R}_+^{2k-1}$is $n$-extendable for every $n$.
\end{proposition}
For a tree, $\tau_n=(\mathcal{V}(\tau_n),\mathcal{E}(\tau_n))$, with $n$ vertices and $n_l$ number of leaves, we consider its $L$-tree, $L(\tau_n,l_1,\ldots,l_k)$, constructed from $k$ leaves, chosen uniformly, with $k <n_l$. The question of extendability is a question of whether models specified in terms of joint distributions over a class of index sets are projective. As in the proof above, consider a probability kernel $p_k$ from $\mathcal{L}_1 \times \cdots \times \mathcal{L}_{k-1}$ to $\mathcal{L}_k$. Noting that $f_n(\tau_n)=f_n(\tau_n,l_1,\ldots,l_{n_l})$, the proof follows along the lines that are identical to those Theorem \ref{projective}.
\qed
\\
%%%%%%%%%%%%%%%%%%%%%%%%%%%%%%%%%
\textbf{Proof of Theorem 4 in paper}.  Let $H_n$ be the Dyck path that corresponds to $\tau_n$. Then, $d(\text{root},v)$ is distributed as $H_n(2nv)$. Since for $0 \leq s \leq 1$, $n^{-1/2}H_n(2ns)$ converges weakly in $C[0,1]$ to $B^{ex}(s)$, we can claim that $n^{-1/2}d(\text{root},v) \overset{d} \to B^{ex}(v)$. Note that this conditional density of the Brownian excursion $\frac{2}{\sigma}B^{ex}$ at time $t \in (0,1)$ is given by (see \cite{takacs})
\begin{equation}\label{excursion_density}
f(t,x)=
\frac{x^2 \sigma^3}{4\sqrt{2\pi t^3(1-t)^3}}e^{\frac{-x^2\sigma^2}{8t(1-t)}}, \quad x>0 .
\end{equation}
Suppose $V$ is uniform on $[0,1]$, and we wish to determine the density of $\frac{2}{\sigma}B^{ex}(V)$, the unconditional density of $B^{ex}(V)$ as
\[
r(x)=\displaystyle \int_0^\infty f(s,x) ds,
\]
since $V$ is uniform on $[0,1]$. Note that the map $B^{ex} \mapsto B^{ex}(V)$ is a one-dimensional random coordinate projection, and is clearly continuous on $C[0,1]$ with respect to the uniform norm. Using the continuous mapping theorem (see \cite{bill}), $d(\text{root},V) \overset{d} \to \frac{2}{\sigma}B^{ex}(V)$, which follows a Rayleigh distribution.
\qed

\section{Computing notes}
We provide details on the parallel and high-performance
simulation platform that we used for our experiments; this software has been
open-sourced on GitHub (\texttt{www.github.com/pkambadu/DyckPaths}) under a
BSD-style license. Our implementation is written in C++ and makes use of the Boost Graph
library (\cite{schling2011boost}) to represent trees, the Boost options library
(\cite{schling2011boost}) to parse command line options, the Boost random
library (\cite{schling2011boost}) to generate various distributions, and OpenMP
for parallelism and therefore, its dependencies. Our code can be compiled and run on any operating system that has a C++
compiler (with or without OpenMP support) as long as the above mentioned Boost
libraries have also been installed; we have tested our implementation on
Darwin 10.7 using GCC 4.2.1 and Ubuntu Linux 2.6.31-23-server using GCC 4.4.1.

\subsection{Generating conditioned Galton--Watson trees with given offspring distribution}
 In order to generate a conditioned Galton--Watson tree $\tau_n$ with offspring distribution $\pi_k$ based on the algorithm in \cite{LD}, it is required to generate a vector  $\Xi=(\xi_1,\ldots,\xi_n)$ where $\xi_i$ are independent copies from $\pi$; subsequently, we are required to rotate $\Xi$ to ensure that $\sum_{i=1}^n \xi_i=n-1$. We shall describe the construction of the conditioned Galton--Watson tree with unit edge lengths. Such a setup implies that (see \cite{LD}):
\begin{enumerate}
\item $\Xi$ is a multinomial random vector with success probabilities determined by $\pi_k$;
\item Elements of $\Xi$ are bounded between $0$ and $n-1$.
\end{enumerate}

We omit details of the algorithm and refer the details to \cite{LD}. However, we provide here a C++ pseudo-code for the generation of the $\Xi$ and consequently the conditioned Galton--Watson. In the code {\small \verb pi_k }  is $\pi_k$; {\small \verb Xi_vec } is  $\Xi$, number of vertices if set to $N$, and \verb xi  is $\xi_i$.

\begin{figure}[!ht]
\centering
\includegraphics[scale=0.5]{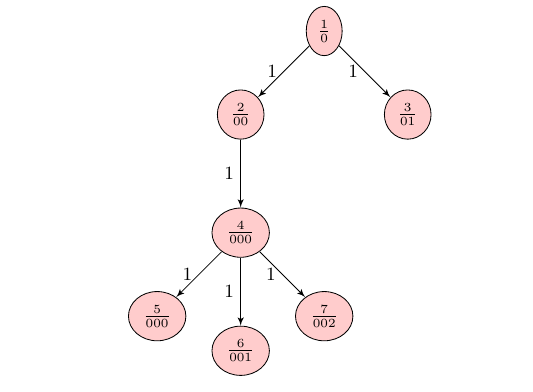}
\caption{
The graph constructed for $\psi = \left[\frac{1}{2}, \frac{2}{1}, \frac{3}{0},
\frac{4}{3}, \frac{5}{0}, \frac{6}{0}, \frac{7}{0}\right]$ using DFS-visit with the DFS-label for each node starting from the root ($0$);
% the
%left-most child  appends a $0$ to it's parent's DFS-label to form it's own
%label, it's right sibling appends a $1$, and so on.
}
\label{fig:construct}
\normalsize
\end{figure}

{\small
\begin{verbatim}
template <typename P_k_type>
std::vector<int> generate_tree (int N, int seed, P_k_type pi_k) {
  std::vector<int> Xi_vec;
  boost::mt19937 prng(seed);

  bool generated(false);
  do {
    int k(0), nodes_consumed(0), edges_consumed(0);
    double sum_of_xi(0.0);
    while (N != nodes_consumed) {
      const double xi = pi_k(i);
      boost::binomial_distribution<> dist((n-sum_of_N_i), xi/(1-sum_of_xi));
      const int n_k = dist(prng);
      for (int i=nodes_consumed; i<(nodes_consumed+n_k);++i) Xi_vec[i]=k;
      nodes_consumed += n_k;
      edges_consumed += n_k*k;
      sum_of_xi += xi;
      ++k;
    }
    if (edges_consumed==(N-1)) generated=true;
  } while (false==generated);

  return Xi_vec;
}
\end{verbatim}
}
\noindent Once the vector $\Xi$ has been generated, it is then required to shuffle it to ensure that $\sum_{i=1}^n \xi_i=n-1$. The first $n_0$ entries of $\Xi$ contain $0$, the next $n_1$ entries contain $1$'s, and so on. We first impart random structure to the conditioned Galton--Watson tree represented by $\Xi$ by a random shuffling or permutation of $\Xi$. We then need to rotate $\Xi$ to ensure that a Depth First Search (DFS) traversal
 will cover all the $n$ nodes. As an example, suppose following the shuffling we are left with $\Xi=\left[0,0,1,2\right]$.
Our DFS based construction algorithm would assign $0 (\psi[0])$ children
to the root node, thereby terminating the tree generation. For this $\Xi$ to be valid for our tree construction, we have to rotate to
get  $\Xi=\left[1,2,0,0\right]$. The index $i, 1\leq i\leq n$ at which $\Xi$ has to be rotated is given
in~\cite{LD}; we give the C++ pseudo-code below:

{\small
\begin{verbatim}
void shuffle_and_rotate (std::vector<int>& Xi_vec) {

  std::random_shuffle (psi_vec.begin(), psi_vec.end());

  size_t N = Xi_vec.size();
  std::vector<int> S(N);

  for (size_t i=0; i<N; ++i) S[i] = ((0==i)?1:S[i-1])+(Xi_vec[i]-1);

  int min_ele=std::numeric_limits<int>::max();
  int min_index=-1;
  for (size_t i=0; i<N; ++i)
    if (min_ele>S[i]) { min_ele=S[i]; min_index=i; }

  std::rotate(Xi_vec.begin(), Xi_vec.begin()+min_index+1, Xi_vec.end());
}
\end{verbatim}
}
Given a properly constructed, shuffled and rotated $\Xi$, construction of the
conditioned Galton--Watson tree is achieved by a DFS based algorithm that is best
illustrated through the use of an example. Consider $\Xi=[2,1,0,3,0,0,0]$; when augmented with the index information,
$\psi = \left[\frac{1}{2}, \frac{2}{1}, \frac{3}{0}, \frac{4}{3}, \frac{5}{0},
\frac{6}{0}, \frac{7}{0}\right]$; here, the numerator denotes the node-ID and the denominator denotes the
number of children (out-degree) of the node. We start by considering node $1$ as the root of the conditioned Galton--Watson tree; in our example,
node 1 has an out-degree of $2$. Therefore, we mark nodes $2$ and $3$ as the children of $1$ and connect them
in our tree. As we explore in DFS-order, we next consider node $2$, which has $1$ child; as
the next unmarked node is $4$, we connect $4$ to be $2$'s child. Next, we explore $4$, which has $3$ children; therefore, we allocate $5,6,7$ as $4$'s children and connect them. Next we explore nodes $5,6,$ and $7$, each of which has $0$ children before returning to node $3$, which also has $0$ children. This completes our tree construction, which is shown in
Figure~\ref{fig:construct}.
 Notice that in addition to constructing the tree, a DFS-label is generated for
every node starting from the root ($0$); the left-most child  appends a $0$ to
it's parent's DFS-label to form it's own label, it's right sibling appends a
$1$, and so on. This DFS-label is useful in finding the LCA of two nodes in a conditioned Galton--Watson tree.
\subsection{Constructing $L$-trees}
Let $T=(V,E)$ represent our conditioned Galton--Watson tree, with vertex set $V$ and
edge set $E$. Let $B$ represent a subset of vertices of $V$ that includes the root.
Given $(T,B)$, the goal of the LCA algorithm is to construct a minimum path
tree $T_{new}$ that contains all the vertices of $B$ and some additional LCA
vertices; that is, we augment $B$ with additional LCA vertices.
Notice that in the worst case, after running our LCA algorithm, $B=V$.
The LCA of two nodes $(v_1,v_2)$ in a tree is a node that is the lowest among
the common ancestors of $v_1$ and $v_2$; the caveat is that the LCA of two
nodes can be one of the nodes themselves as each node is it's own ancestor.
The algorithm to compute the LCA is rather simple:
\begin{enumerate}
\item Let $V_{new}$ represent the list of vertices $B$ plus the LCA vertices that
augment $B$ --- initialize this list to $V_{new}=B$;
\item
For each pair of vertices $(v_1,v_2)\in B$, compute $v_{LCA}$, the LCA of
$(v_1,v_2)$ and add it to $V_{new}$; there are $ \binom {\#B}{2}$ such vertex
pairs;
\item
Construct the $L$-tree $T_{new}$ by joining the vertices in
$V_{new}$ using the edge information in $T$; specifically, when connecting
vertices that originally had an intermediate vertex between them in $T$,
augment the new edge to include the weights of the edges that were skipped in
$T$.
\end{enumerate}
We now turn our attention to efficient computation of $v_{LCA}=L(v_1,v_2)$. Notice that we label each of the vertices in $T$ with their DFS-label (see
Figure~\ref{fig:construct}).
This DFS-label can be used directly to determine the LCA; the LCA of
$(v_1,v_2)$ is the longest common prefix of the labels of $v_1$ and $v_2$.
For example, consider the nodes $3$ and $5$ in Figure~\ref{fig:construct},
which have the labels $`01'$ and $`000'$, respectively.
The longest common prefix is $`0'$, which points to vertex $0$, which also is
the LCA of $3$ and $5$.
As we store the DFS-label of each node succinctly as a string, we are able to
quickly find the LCA using the \verb std::mismatch  algorithm, which returns
the first position of mismatch in the two DFS-labels.\\

\subsection{Parallel execution}
The basic control structure of our simulations is: (a) generate a large number
of conditioned Galton--Watson trees; (b) compute local statistics on each conditioned Galton--Watson tree; and (c) combine
the local statistics to make inferences.
As mentioned earlier, generating a single tree is expensive and may potentially
incur many failed attempts before success.
Therefore, we parallelize the simulation framework by parallelizing step (b)
above using OpenMP; that is, multiple trials of the experiments are run
simultaneously when possible and combined with care to ensure consistency.
Given that most of the computing hardware has inherent parallelism in the form
of multi-cores and multi-sockets, our approach results in linear speedups
(w.r.t number of computational resources) in throughput.
Notice that parallelizing step (a) is hard both because of the sequential
dependency in generating $\Xi$ from the multinomial distribution and because
our current random number generators are not thread-safe.
However, as we conduct thousands of experiments, we are able to fully utilize
clusters with similar processor counts; that is, parallelizing step (a) is not
necessary.

\section{Simulations}
\subsection{Poisson process model}
We simulate binary trees from the non-homogeneous Poisson process model with intensity $\lambda(t)=2.1t$, for $t>0$. We refer to this model as the model under $H_0$. Similarly, trees with intensity $\lambda(t)=4.9t$ are generated, and rejection rates for the tests in Theorem 1 are tabulated in Table \ref{binary_test}. In each test, rejection rates were computed by averaging over multiple permutations of the leaves (see p. 219, Theorem 6.3.1 of \citet{LR}) , verifying the permutational invariance of the test. 
\begin{table}[!htb]
\small
\begin{tabular}{|c|cc|cc|cc||cc|cc|cc|}
%\hline
\multicolumn{7}{c}{One-sample}&\multicolumn{6}{c}{Two-sample}\\
\hline
\multirow{2}{*}{Tree sizes}& \multicolumn{2}{c|}{$N=10$}&\multicolumn{2}{c|}{$N=20$}&\multicolumn{2}{c||}{$N=50$}& 
\multicolumn{2}{c|}{$N=10$}&\multicolumn{2}{c|}{$N=20$}&\multicolumn{2}{c|}{$N=50$}\\
%\hline
&$H_0$& $H_a$&$H_0$& $H_a$&$H_0$& $H_a$&$H_0$& $H_a$&$H_0$& $H_a$&$H_0$& $H_a$\\
\hline
10&0.09&0.79&0.05&0.88&0.03&0.95& 0.09&0.81&0.06&0.91&0.04&0.94\\
20&0.06&0.93&0.04&0.93&0.05&0.95& 0.07&0.89&0.05&0.89&0.01&0.95\\
50&0.07&0.91&0.04&0.95&0.03&0.94&0.04&0.94&0.06&0.94&0.05&0.93\\
100&0.05&0.94&0.05&0.96&0.04&0.95&0.05&0.93&0.04&0.92&0.055&0.95\\
\hline
\end{tabular}
\vspace*{4mm}
%\end{table}
%\vspace*{4mm}
\caption{\footnotesize{Tests for binary trees from Poisson process model with intensity $\lambda(t)=\theta t$. Here $H_0:\theta=2.1$ and $H_a:\theta=4.9$. Rejection rates, at $\alpha=0.05$, computed by averaging over multiple permutation of the leaves, are provided under $H_0$ and $H_a$ for varying sample sizes $N$ and tree sizes. }}
\label{binary_test}
\end{table}
\subsection{Galton--Watson trees}
For each tree, the branch lengths of the conditioned Galton--Watson trees are generated from a uniform distribution between 0 and 2. Although the framework allows each tree to contain a different number of vertices, for simplicity, we ensure that each generated tree is of the same size. However, we examine the performance of the tests at different sample sizes $N$ since we are interested in the sample sizes for consistent estimation of $\sigma^2$. As the proposed framework is based on trees with a large number of vertices, following some trials, we find that the performance of the tests is reasonable, with trees containing around 1000 vertices. Therefore, all simulation exercises for conditioned Galton--Watson trees are carried out with trees that contain 1000 vertices.

Next, we generate conditioned Galton--Watson trees with varying sample sizes from some of the distributions listed in Proposition 1 of the main document, with each containing 1000 vertices. Each simulation exercise contains 1000 simulated trials. We first compute the estimator of $\sigma^2$ from Theorem 4 of the main document based on choosing a vertex at random; then, we construct $L$-trees by randomly choosing 25 leaves in each tree. We compute the total path length of the $L$-trees and then compute the test statistic by scaling the total path length by the estimator. Rejection rates are computed by averaging over multiple permutations of the 25 leaves, verifying the invariance property of the test statistic. We compare the performance of the $\chi^2$ and $F$ tests to that of Permutation tests, with 5000 permutations, based on the same test statistic used in Table 1 in the main paper.

In order to examine rejection rates under the alternative hypotheses, we generate unconditioned Galton--Watson trees, and two types of phylogenetic trees that do not belong to the class of conditioned Galton--Watson trees: trees based on a birth-death process with speciation rate of 2 on fixed taxa, and genealogical trees of Kingman's coalescent processes. However, there is a subtle connection between the coalescent process and the CRT  \citep{CH}, which manifests itself through a poor rejection rate in the tests. 

\begin{figure}[!ht]
\begin{center}
\begin{tabular}{c c c }
\includegraphics[scale=0.25]{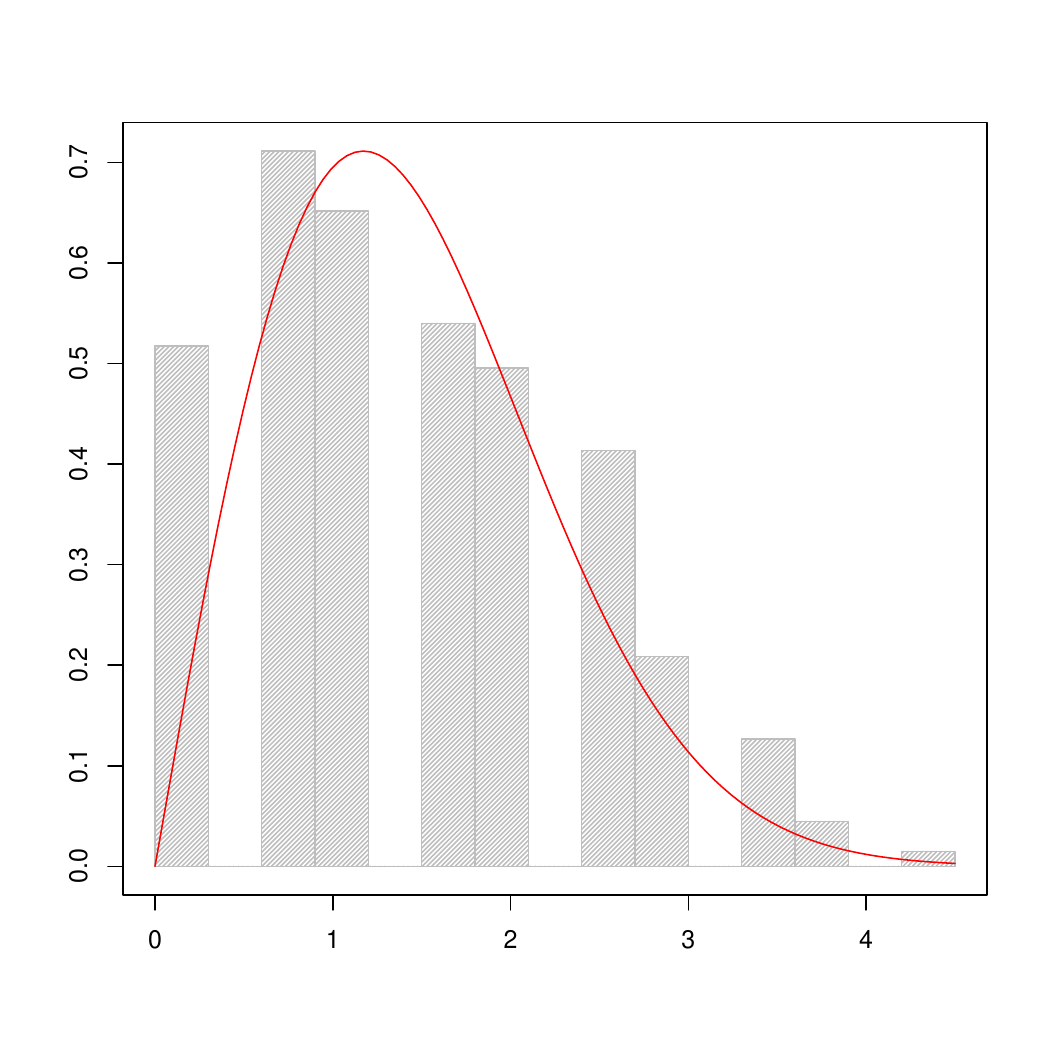} & \includegraphics[scale=0.25]{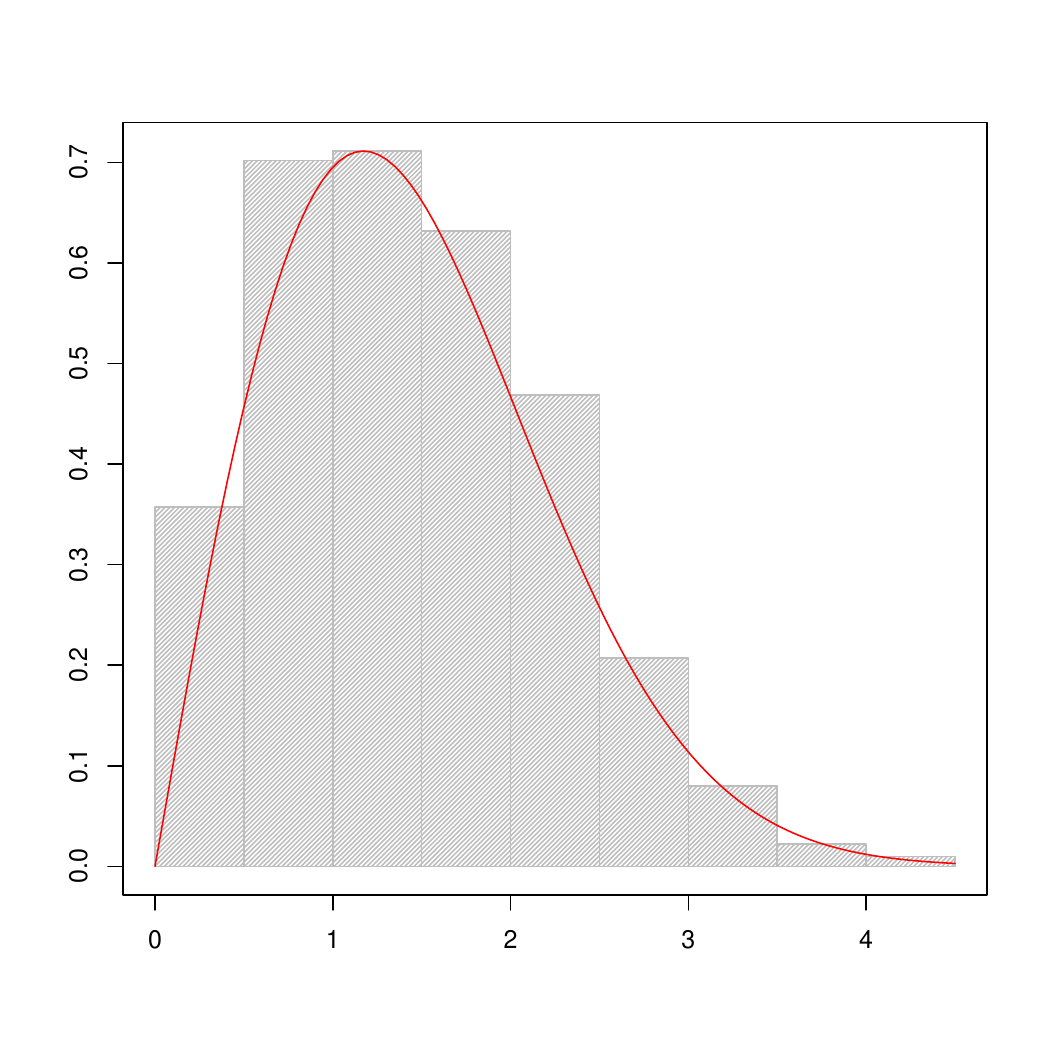}
&\includegraphics[scale=0.25]{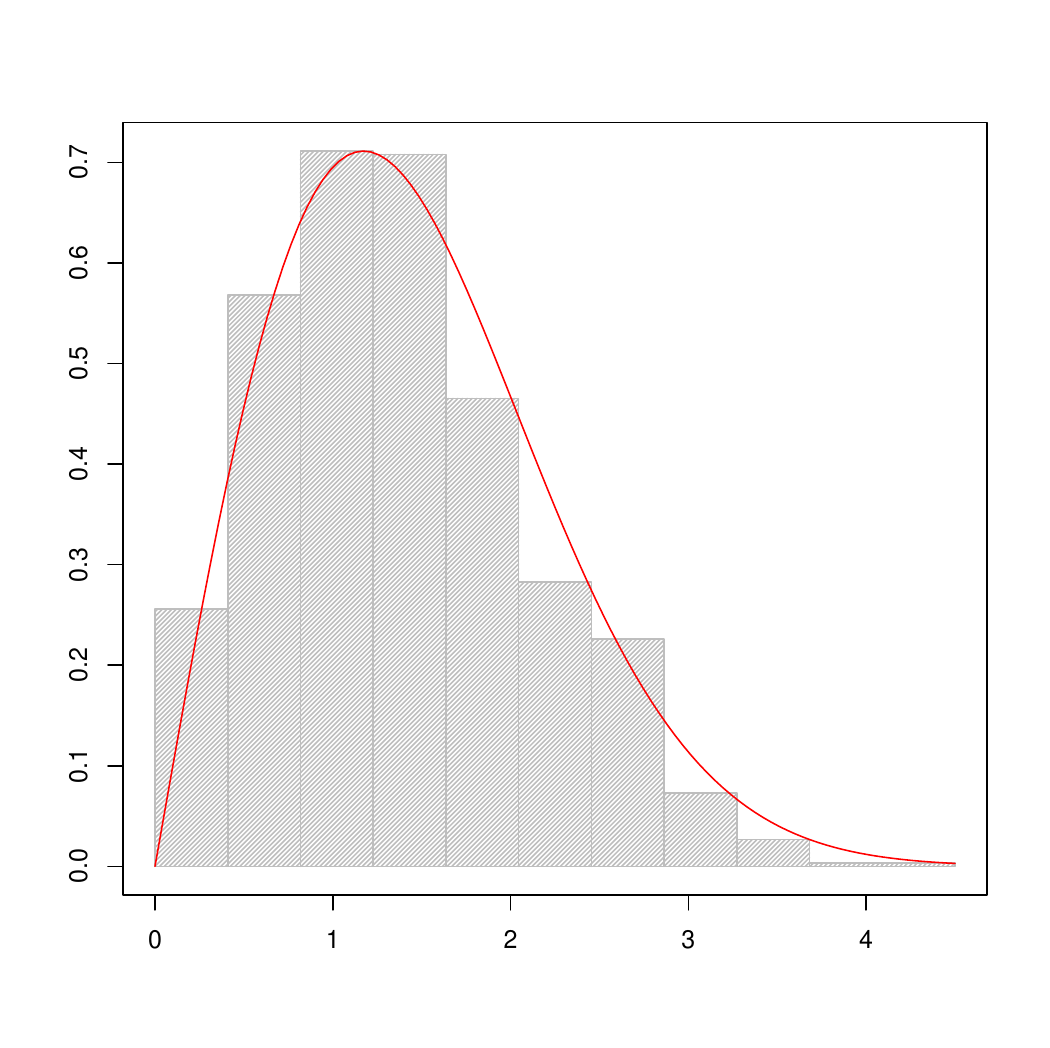}
\end{tabular}
\end{center}
\caption{\footnotesize{Histograms of $n^{-1/2}d(\text{root},V)$ where $V$ is vertex chosen at random on 10000 CGW trees from a Bin(2,1/2) distribution with respective number of vertices $n=$ 10, 100 and 1000 (from left). Solid red curve is the Rayleigh density with scale $\sqrt{2}$.}}
\label{histheight}
\end{figure}

Asymptotic verification of the results from Theorem 4 in the main document is carried out by simulating conditioned Galton--Watson trees with Bin(2,1/2) offspring distribution and $\sigma^2=1/2$. Figure \ref{histheight} plots histograms of the normalized distance of a randomly chosen vertex from trees with varying number of vertices.  
\begin{figure}[!ht]
	\footnotesize
	\parbox{0.45\linewidth}{
		\begin{tabular}{c c }
			\includegraphics[scale=0.3]{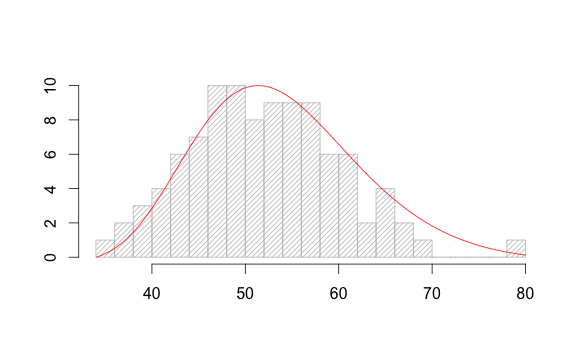} & \includegraphics[scale=0.3]{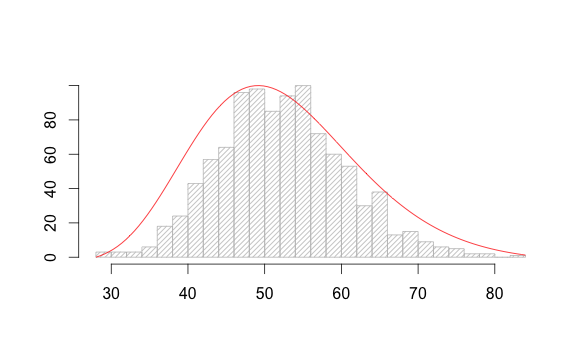}
			%&\includegraphics[scale=0.25]{figs/height_geo1000}
		\end{tabular}
	}
	\hfill
	\parbox{0.25\linewidth}{
		\begin{tabular}{|c c|}
			\hline
			$N$ & $\hat{\sigma}^2$ \\
			\hline
			50& 1.94\\
			100& 1.96\\
			500& 1.98\\
			1000& 2.01\\
			\hline
		\end{tabular}
	}
	\caption{\footnotesize Histograms of total path length, scaled by the estimate of $\sigma^2$, of $L$-trees of 100 (left histogram) and 1000 (right histogram) conditioned Galton--Watson trees from Geo$(1/2)$ with 1000 vertices each. $L$-trees were constructed by randomly choosing 25 leaves. Solid red curve is the Gamma density with scale parameter 2 and shape 25. The box lists the accuracy of the estimator of the variance of Geo(1/2), $\sigma^2=2$, at different sample sizes.}
	\label{histLCA}
\end{figure}
Recall from Theorem 2 of the main document that the critical function for the test is based on the statistic that corresponds to the scaled total path length or sum of branch lengths of the $L$-tree constructed from a randomly chosen subset of the leaves of the conditioned Galton--Watson tree; the scaling constant was estimated using the estimator from Theorem 4 of the main document. Figure \ref{histLCA} compares the histogram of the scaled total path length of the constructed $L$-tree and the Gamma distribution from Proposition \ref{gamma}, and also reports the accuracy of the estimator; we see that the approximation is reasonable for trees with a large number of vertices.
\end{document}